\documentclass[12pt]{amsart}
\usepackage[margin=3.3cm]{geometry}
\usepackage{amssymb,epsfig}
\usepackage[all]{xy}

\newcommand{\Z}{\mathbb{Z}}
\newcommand{\A}{\mathcal{A}}

\newcommand{\D}{\mathcal{D}}
\renewcommand{\H}{\mathcal{H}}

\newcommand{\sqbinom}[2]{\begin{bmatrix} #1\\#2 \end{bmatrix}}

\renewcommand{\a}{{\alpha}}
\newcommand{\be}{{\beta}}
\newcommand{\g}{{\gamma}}

\renewcommand{\d}{{\delta}}
\newcommand{\De}{{\Delta}}
\newcommand{\e}{{\varepsilon}}
\newcommand{\z}{{\zeta}}

\renewcommand{\l}{{\lambda}}
\newcommand{\m}{{\mu}}
\newcommand{\s}{{\sigma}}
\renewcommand{\t}{{\tau}}
\newcommand{\w}{{\omega}}

\newcommand{\W}{{\widehat W}}

\newcommand{\ra}{{\ \longrightarrow \ }}

\newcommand{\aut}{\operatorname{Aut}}

\newcommand{\qv}{{\mathbb{Q}(v)}}

\newcommand{\eugn}{{\widehat{\mathbf{U}}(\widehat{\mathfrak{gl}_n})}}

\newcommand{\ugn}{{\mathbf{U}(\widehat{\mathfrak{gl}_n})}}

\newcommand{\ugln}{\mathbf{U}(\mathfrak{gl}_n)}

\newcommand{\gl}{\mathfrak{gl}}
\newcommand{\sqnr}{{\widehat{\mathbf{S}}_q(n, r)}}
\newcommand{\svnr}{{\widehat{\mathbf{S}}_v(n, r)}}

\newcommand{\zed}{\mathbb{Z}}

\renewcommand{\Im}{\operatorname{Im}}
\newcommand{\End}{\operatorname{End}}
\newcommand{\Hom}{\operatorname{Hom}}

\newcommand{\ii}{\mathtt{i}}
\newcommand{\EE}{\mathsf{E}}
\newcommand{\FF}{\mathsf{F}}
\newcommand{\KK}{\mathsf{K}}
\newcommand{\iii}{\mathtt{i}}
\newcommand{\eee}{\mathsf{e}}
\newcommand{\fff}{\mathsf{f}}
\newcommand{\shatnr}{\widehat{S}(n,r)}

\newcommand{\bS}{\mathbf{S}}

\newcommand{\seq}{\mathfrak{S}_{r,n}}
\newcommand{\mat}{\mathfrak{S}_{r,n,n}}
\newcommand{\matap}{\mathfrak{S}^{\text{ap}}_{r,n,n}}

\newtheorem{thm}[figure]{Theorem}
\newtheorem{lem}[figure]{Lemma}
\newtheorem{prop}[figure]{Proposition}
\newtheorem{cor}[figure]{Corollary}
\theoremstyle{definition}
\newtheorem{defn}[figure]{Definition}
\newtheorem{rmk}[figure]{Remark}

\newtheorem{exam}[figure]{Example}
\numberwithin{figure}{subsection}

\allowdisplaybreaks

\begin{document}
\title{Presenting affine $q$-Schur algebras}

\author{S.~R.~Doty}
\address{Department of Mathematics and Statistics,
Loyola University Chicago,
Chicago, IL 60626 U.S.A.}
\email{sdoty@luc.edu}
\thanks{}

\author{R.~M.~Green} 
\address{Department of Mathematics, University of
Colorado, Campus Box 395, Boulder, CO 80309 U.S.A.}
\email{rmg@euclid.colorado.edu}
\thanks{}

\begin{abstract}  
We obtain a presentation of certain affine $q$-Schur algebras in terms
of generators and relations.  The presentation is obtained by adding
more relations to the usual presentation of the quantized enveloping
algebra of type affine $\gl_n$.  Our results extend and rely on the
corresponding result for the $q$-Schur algebra of the symmetric group,
which were proved by the first author and Giaquinto.
\end{abstract}

\maketitle

\section*{Introduction}

Let $V'$ be a vector space of finite dimension $n$. On the tensor
space $(V')^{\otimes r}$ we have natural commuting actions of the
general linear group $\mathsf{GL}(V')$ and the symmetric group
$\mathcal{S}_r$.  Schur observed that the centralizer algebra of each
action equals the image of the other action in $\End((V')^{\otimes
r})$, in characteristic zero, and Schur and Weyl used this observation
to transfer information about the representations of $\mathcal{S}_r$
to information about the representations of $\mathsf{GL}(V')$. That
this Schur--Weyl duality holds in arbitrary characteristic was first
observed in \cite{Dec-Pro}, although a special case was already used
in \cite{Car-Lus}.  In recent years, there have appeared various
applications of the Schur--Weyl duality viewpoint to modular
representations. The Schur algebras $S(n,r)$ first defined in
\cite{JAGreen} play a fundamental role in such interactions.

Jimbo \cite{Jimbo} and (independently) Dipper and James \cite{DJ2}
observed that the tensor space $(V')^{\otimes r}$ has a $q$-analogue
in which the mutually centralizing actions of $\mathsf{GL}(V')$ and
$\mathcal{S}_r$ become mutually centralizing actions of a quantized
enveloping algebra $\ugln$ and of the Iwahori-Hecke algebra
$\H(\mathcal{S}_r)$ corresponding to $\mathcal{S}_r$. In this context,
the ordinary Schur algebra $S(n,r)$ is replaced by the $q$-Schur
algebra $\bS_q(n,r)$. Dipper and James also showed that the $q$-Schur
algebras determine the representations of finite general linear groups
in non-defining characteristic.
 
An affine version of Schur--Weyl duality was first described in
\cite{CP}. A different version, in which the vector space $V'$ is
replaced by an infinite dimensional vector space $V$, is given in
\cite{G15}, and we follow the latter approach here.  In the affine
(type $A$) setting, the mutually commuting actions are of an affine
quantized enveloping algebra $\ugn$ and an extended affine Hecke
algebra $\H(\W)$ corresponding to an extended affine Weyl group $\W$
containing the affine Weyl group $W$ of type $\widehat{A}_{r-1}$. The
affine $q$-Schur algebra $\sqnr$ in this context, which is also
infinite dimensional, was first studied in \cite{G15}, \cite{L11}, and
\cite{VV}.

Recently, a new approach to Schur algebras or their $q$-analogues was
given in \cite{DG}, where it was shown that they may be defined by
generators and relations in a manner compatible with the usual
defining presentation of the enveloping algebra or its corresponding
quantized enveloping algebra. The purpose of this paper is to extend
that result to the affine case --- that is, to describe the affine
$q$-Schur algebra $\sqnr$ by generators and relations compatible with
the defining presentation of $\ugn$. This result is formulated in
Theorem \ref{Thm:1.6.1}, under the assumption that $n>r$.  An
equivalent result, which describes the affine $q$-Schur algebra as a
quotient of Lusztig's modified form of the quantized enveloping
algebra, is given in Theorem \ref{Thm:2.6.1}. These results depend on
a different presentation, also valid for $n>r$, of the $q$-Schur
algebra given in \cite[Proposition 2.5.1]{G15}. A different approach
to the results of this paper seems to be indicated for the case $n \le
r$.

The organization of the paper is as follows. In Section 1 we give
necessary background information, and formulate our main result. In
Section 2 we give the proof of Theorem \ref{Thm:1.6.1}, and we also
give, in Section \ref{ss:2.6}, the alternative presentation mentioned
above. Finally, in Section 3 we outline the analogous results in the
classical case, when the quantum parameter is specialized to 1.

\section{Preliminaries and statement of main results}

Our main result, stated in \S\ref{ss:1.6}, is a presentation by
generators and relations of the affine $q$-Schur algebra.  In order to
put this result in context, we review some of the definitions of the
algebra that have been given in the literature.

\subsection{Affine Weyl groups of type $A$}\label{ss:1.1}

The affine Weyl group will play a key role, both in our definitions
and our methods of proof, so we define it first.

\begin{figure}[ht]  
\centering
\includegraphics{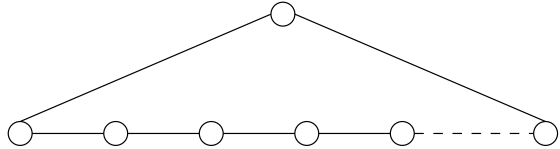}
\caption{Dynkin diagram of type $\widehat{A}_{r-1}$}
\label{fig:affa}       
\end{figure}

The Weyl group we consider in this paper is that of type
$\widehat{A}_{r-1}$, where we intend $r \ge 3$.  This corresponds to
the Dynkin diagram in Figure \ref{fig:affa}.

The number of vertices in the graph in Figure \ref{fig:affa} is $r$,
as the top vertex (numbered $r$) is regarded as an extra relative to
the remainder of the graph, which is a Coxeter graph of type
$A_{r-1}$.

We associate a Weyl group, $W = W(\widehat{A}_{r-1})$, to this Dynkin
diagram in the usual way (as in \cite[\S2.1]{H}).  This associates to
node $i$ of the graph a generating involution $s_i$ of $W$, where $s_i
s_j = s_j s_i$ if $i$ and $j$ are not connected in the graph, and 
$$
s_i s_j s_i = s_j s_i s_j
$$ 
if $i$ and $j$ are connected in the graph.  For $t \in \Z$, it is
convenient to denote by $\overline{t}$ the congruence class of $t$
modulo $r$, taking values in the set $\{1, 2, \ldots, r\}$.  For the
purposes of this paper, it is helpful to think of the group $W$ as
follows, based on a result of Lusztig \cite{L6}. (Note that we write
maps on the right when dealing with permutations.)

\begin{prop}\label{Prop:1.1.2}
There exists a group isomorphism from $W$ to the set of permutations
of $\Z$ satisfying the following conditions:
\begin{align}
(i+r)w &= (i)w + r  \tag{a}\\
\sum_{t = 1}^r (t)w &= \sum_{t = 1}^r t  \tag{b}
\end{align}
such that $s_i$ is mapped to the permutation
$$
t \mapsto 
\begin{cases}
t & \text{ if } \overline{t} \ne \overline{i}, \overline{i + 1},\\
t-1 & \text{ if } \overline{t} = \overline{i + 1},\\
t+1 & \text{ if } \overline{t} = \overline{i},
\end{cases}
$$ 
for all $t \in \Z$.
\end{prop}

For reasons relating to weight spaces which will become clear later,
we consider a larger group $\W$ of permutations of $\Z$.

\begin{defn}\label{Def:1.1.3}
Let $\rho$ be the permutation of $\Z$ taking $t$ to $t + 1$ for
all $t$.  Then the group $\W$ is defined to be the subgroup of
permutations of $\Z$ generated by the group $W$ and $\rho$.
\end{defn}

As will become clear later, the point of $\rho$ is that conjugation by
$\rho$ will correspond to a graph automorphism of the Dynkin diagram
given by rotation by one place.

\begin{prop}\label{Prop:1.1.4}
{\rm (i)} There exists a group isomorphism from $\W$ to the set of
permutations of $\Z$ satisfying the following conditions:
\begin{align}
(i+r)w &= (i)w + r  \tag{a}\\
\sum_{t = 1}^r (t)w &\equiv \sum_{t = 1}^r t \mod r.\tag{b}
\end{align}
\par\noindent{\rm (ii)} Any element of $\W$ is uniquely expressible in
the form $\rho^z w$ for $z \in \Z$ and $w \in W$.  Conversely, any
element of this form is an element of $\W$.  
\par\noindent{\rm (iii)}
Let $S \cong \mathcal{S}_r$ be the subgroup of $\W$ generated by
$$
\{s_1, s_2, \ldots, s_{r-1}\}.
$$  
Let $Z$ be the subgroup of $\W$ consisting of all permutations 
$z$ satisfying 
$$
(t)z \equiv t \mod r
$$ 
for all $t$.  Then $\Z^r \cong Z \triangleleft \W$ and $\W$
is the semidirect product of $S$ and $Z$.
\end{prop}

\begin{proof}
The three parts are proved in \cite[Proposition 1.1.3, Corollary
1.1.4, Proposition 1.1.5]{G15} respectively.
\end{proof}

It is convenient to extend the usual notion of the length of an
element of a Coxeter group to the group $\W$ in the following way.

\begin{defn}\label{Def:1.1.5}
For $w \in W$ the length $\ell(w)$ of $w$ is the length of a word of
minimal length in the group generators $s_i$ of $W$ which is equal to
$w$.  The length, $\ell(w')$, of a typical element $w' = \rho^z w$ of
$\W$ (where $z \in \Z$ and $w \in W$) is defined to be $\ell(w)$.
\end{defn}

When the affine Weyl group is thought of in the above way, the
familiar notions of length and distinguished coset representatives may
be adapted from the corresponding notions for Coxeter groups.

\begin{defn}\label{Def:1.1.6}
Let $\Pi$ be the set of subsets of $S = \{s_1, s_2, \ldots, s_r\}$,
excluding $S$ itself.  For each $\pi \in \Pi$, we define the subgroup
$\W_\pi$ of $\W$ to be that generated by $\{s_i \in \pi\}$.  (Such a
subgroup is called a parabolic subgroup.)  We will sometimes write
$W_\pi$ for $\W_\pi$ to emphasize that it is a subgroup of $W$.  Let
$\Pi'$ be the set of elements of $\Pi$ that omit the generator $s_r$.
\end{defn}

All the subgroups $\W_\pi$ are subgroups of $W$, and are parabolic
subgroups in the usual sense of Coxeter groups.  Furthermore, each
such $\W_\pi$ is isomorphic to a direct product of Coxeter groups of
type $A$ (i.e., finite symmetric groups) corresponding to the
connected components of the Dynkin diagram obtained after omitting the
elements $s_i$ which do not occur in $\pi$.  We will appeal to these
facts freely in the sequel.

\begin{defn}\label{Def:1.1.7}
Let $\pi \in \Pi$.  The subset $\D_\pi$ of $\W$ is the set of
those elements such that for any $w \in \W_\pi$ and $d \in \D_\pi$, 
$$
\ell(w d) = \ell(w) + \ell(d).
$$ 
We call $\D_\pi$ the set of distinguished right coset
representatives of $\W_\pi$ in $\W$.

The subset $\D_\pi^{-1}$ is called the set of distinguished left coset
representatives of $\W_\pi$ in $\W$; elements $d \in \D_\pi^{-1}$ have
the property that $\ell(d w) = \ell(d) + \ell(w)$ for any $w \in
\W_\pi$.
\end{defn}

\begin{prop}\label{Prop:1.1.8}
{\rm (i)} Let $\pi \in \Pi$ and $w \in \W$.  Then $w = w_\pi w^\pi$
for a unique $w_\pi \in \W_\pi$ and $w^\pi \in \D_\pi$.
\par\noindent{\rm (ii)} Let $\pi' \in \Pi$ and $w \in \W$.  Then $w =
w^{\pi'} w_{\pi'}$ for a unique $w_{\pi'} \in \W_{\pi'}$ and $w^{\pi'}
\in \D_{\pi'}$.  
\par\noindent{\rm (iii)} Let $\pi_1, \pi_2 \in \Pi$.
The set $\D_{\pi_1, \pi_2} := \D_{\pi_1} \cap \D_{\pi_2}^{-1}$ is an
irredundantly described set of double $\W_{\pi_1}$--$\W_{\pi_2}$-coset
representatives, each of minimal length in its double coset.
\end{prop}

\begin{proof}
See \cite[Propositions 1.4.4, 1.4.5]{G15}.
\end{proof}

\subsection{Affine Hecke algebras of type $A$}\label{ss:1.2}

We now define the extended affine Hecke algebra $\H = \H(\W)$ of type
$A$.  The Hecke algebra is a $q$-analogue of the group algebra of
$\W$, and is related to $\W$ in the same way as the Hecke algebra
$\H(\mathcal{S}_r)$ of type $A$ is related to the symmetric group
$\mathcal{S}_r$.  In particular, one can recover the group algebra of
$\W$ by replacing the parameter $q$ occurring in the definition of
$\H(\W)$ by $1$.

\begin{defn}\label{Def:1.2.1}
The affine Hecke algebra $\H = \H(\W)$ over $\Z[q, q^{-1}]$ is
the associative, unital algebra with algebra generators 
$$
\{T_{s_1}, \ldots, T_{s_r}\} \cup \{ T_\rho, T_\rho^{-1} \}
$$ 
and relations 
\begin{align}
&T_s^2 = q T_s + (q-1),  \tag{1}\\
&T_s T_t = T_t T_s \text{ if $s$ and $t$ are not adjacent in the
Dynkin diagram},  \tag{2}\\
& T_s T_t T_s = T_t T_s T_t \text{ if $s$ and $t$ are adjacent in the
Dynkin diagram},  \tag{3}\\
&T_\rho T_{s_{i+1}} T_\rho^{-1} = T_{s_i}.  \tag{4}
\end{align}  
In relation (4), we interpret $s_{r+1}$ to mean $s_1$.
\end{defn}

The algebra $\H$ has a better known presentation, known as the
Bernstein presentation, but this is not convenient for our purposes.
The equivalence of the two presentations is well known, and a proof
may be found, for example, in \cite[Theorem 4.2.5]{G15}.  However, it
will be convenient to have the following modified version of the
presentation in Definition \ref{Def:1.2.1}.

\begin{lem}\label{Lem:1.2.2}
The affine Hecke algebra $\H = \H(\W)$ over $\Z[q, q^{-1}]$ is
the associative, unital algebra with algebra generators 
$$
\{T_{s_1}, \ldots, T_{s_{r-1}}\} \cup \{ T_\rho, T_\rho^{-1} \}
$$ 
and relations 
\begin{align}
&T_{s_i}^2 = q T_{s_i} + (q-1),  \tag{$1^\prime$} \\
&T_{s_i} T_{s_j} = T_{s_j} T_{s_i} \text{ if $|i - j| > 1$},  
 \tag{$2^\prime$}\\
&T_{s_i} T_{s_j} T_{s_i} = T_{s_j} T_{s_i} T_{s_j} \text{ if $|i - j| = 1$},  
 \tag{$3^\prime$}\\
&T_\rho T_{s_{i+1}} T_\rho^{-1} = T_{s_{i}} \text{ if $1 \le i < r-1$},  
 \tag{$4^\prime$} \\
&T_\rho^r T_{s_{i}} T_\rho^{-r} = T_{s_{i}} \text{ if $1 \le i \le r-1$},
 \tag{$5^\prime$}
\end{align}
\end{lem}

\begin{proof}
It is clear that relations ($1^\prime$)--($5^\prime$) are consequences
of relations (1)--(4).  For the converse direction, we define $T_{s_r}
:= T_\rho T_{s_1} T_{\rho}^{-1}$; the remaining cases of relations
(1)--(4) may then be obtained from relations
($1^\prime$)--($5^\prime$) by conjugating by $T_{\rho}$ or by
$T_{\rho}^{-1}$.
\end{proof}

\begin{defn}\label{Def:1.2.3}
Let $w \in W$.  The element $T_w$ of $\H(W)$ is defined as 
$$
T_{s_{i_1}} \cdots T_{s_{i_m}},
$$ 
where $s_{i_1} \cdots s_{i_m}$ is a reduced expression for $w$
(i.e., one with $m$ minimal).  (This is well-defined by standard
properties of Coxeter groups.)

If $w' \in \W$ is of form $\rho^z w$ for $w \in W$, we denote by
$T_{w'}$ the element $T_\rho^z T_w.$ (This is well-defined by
Proposition \ref{Prop:1.1.4} (ii).)
\end{defn}

\begin{prop}\label{Prop:1.2.4}
{\rm (i)} A free $\Z[q, q^{-1}]$-basis for $\H$ is given by the
set $\{T_w : w \in \W\}.$
\par\noindent{\rm (ii)} As a $\Z[q, q^{-1}]$-algebra, $\H$ is generated by
$T_{s_1}$, $T_\rho$ and $T_\rho^{-1}$.
\end{prop}

\begin{proof}
See \cite[Proposition 1.2.3, Lemma 1.2.4]{G15}.
\end{proof}

\subsection{The affine $q$-Schur algebra as an endomorphism algebra}
\label{ss:1.3}

We first present the definition of the affine $q$-Schur algebra as
given in \cite[\S2]{G15}.

\begin{defn}\label{Def:1.3.1}
A weight is a composition $\lambda = (\lambda_1, \lambda_2, \ldots,
\lambda_n)$ of $r$ into $n$ pieces, that is, a finite sequence of
nonnegative integers whose sum is $r$.  (There is no monotonicity
assumption on the sequence.)  We denote the set of weights by
$\Lambda(n, r)$.

The $r$-tuple $\ell(\lambda)$ of a weight $\lambda$ is the weakly
increasing sequence of integers where there are $\lambda_i$
occurrences of the entry $i$.

The Young subgroup $\mathcal{S}_{\lambda} \subseteq \mathcal{S}_r
\subseteq W \subseteq \W$ is the subgroup of permutations of the set $
\{1, 2, \ldots, r\} $ that leaves invariant the following sets of
integers: 
$$ 
\{ 1, 2, \ldots, \lambda_1 \}, \{\lambda_1 + 1, \lambda_1
+ 2, \ldots, \lambda_1 + \lambda_2\}, \{\lambda_1 + \lambda_2 + 1,
\ldots \}, \ldots .
$$

The weight $\omega$ is given by the $n$-tuple 
\[
( \underbrace{1, 1, \ldots, 1}_r, \underbrace{0, 0, \ldots, 0}_{n-r} ).
\]
\end{defn}

\begin{rmk}\label{Rem:1.3.2}
The Young subgroup $\mathcal{S}_\lambda \subseteq \mathcal{S}_r$ can
be thought of as a group $\W_\lambda$ for some $\lambda \in \Pi'$.
Note, however, that different compositions $\lambda$ can give rise to
canonically isomorphic groups.  Also note that we require $n \ge r$
for $\omega$ to exist.
\end{rmk}

\begin{defn}\label{Def:1.3.3}
Let $\lambda \in \Pi$.  
For $t \in \Z$, the parabolic subgroup $\W_{\lambda + t}$ is
the one generated by those elements $s_{\overline{i + t}}$ where $i$
is such that $s_i$ lies in $\W_{\lambda}$.
We also use the notation $\D_{\lambda + t}$ with the obvious meaning.

The element $x_{\lambda + t} \in \H$ is defined as $$
x_{\lambda + t} := \sum_{w \in \W_{\lambda + t}} T_w
.$$  We will write $x_\lambda$ for $x_{\lambda + 0}$.
\end{defn}

\begin{defn}\label{Def:1.3.4}
The affine $q$-Schur algebra $\sqnr$ over $\Z[q, q^{-1}]$ 
is defined by $$
\sqnr := \End_\H \left( \bigoplus_{\l \in \Lambda(n, r)} x_\l \H \right),
$$ where $\H = \H(\W)$.
\end{defn}

There is a basis for $\sqnr$ similar to Dipper and James' basis for the
ordinary $q$-Schur algebra.

\begin{defn}\label{Def:1.3.5}
Let $d \in \W$ be an element of $\mathcal{D}_{\l, \mu}$.  Write $d =
\rho^z c$ (as in Proposition \ref{Prop:1.1.4} (ii)) with $c \in W$.
Then the element
$$\phi_{\l, \mu}^d \in \Hom(x_\mu \H(\W), x_\l \H(\W))$$ 
is defined as 
\begin{align*}
\phi_{\l, \mu}^d (x_\mu) :=& 
\sum_{d' \in \mathcal{D}_\nu \cap W_\mu} x_\l T_\rho^z T_{c d'}\\
=& \sum_{d' \in \mathcal{D}_\nu \cap W_\mu} T_\rho^z x_{\l + z} T_{c d'} = 
\sum_{w \in W_{\l + z} c W_\mu} T_\rho^z T_w
= \sum_{w \in W_\l d W_\mu} T_w
\end{align*}
where $\nu$ is the composition of $n$ corresponding to the standard
Young subgroup $$d^{-1} W_{\l + z} d \cap W_\mu$$ of $W$.
\end{defn}

\begin{thm}\label{Thm:1.3.6}
{\rm (i)} A free $\Z[q, q^{-1}]$-basis for $\sqnr$ is given by
the set 
$$ 
\{\phi_{\l, \mu}^d : \l, \mu \in \Lambda(n, r), \ d \in
\mathcal{D}_{\l, \mu}\} .$$ 
\par\noindent{\rm (ii)} The set of
basis elements $$ \{\phi_{\l, \mu}^d : \l, \mu \in \Lambda(n, r), d
\in \mathcal{S}_r \cap \D_{\l, \mu}\}
$$ 
spans a subalgebra of $\sqnr$ canonically isomorphic to the
$q$-Schur algebra $\bS_q(n, r)$.
\par\noindent{\rm (iii)}
The set of basis elements $$
\{\phi_{\w, \w}^d : d \in \W\}
$$ spans a subalgebra canonically isomorphic to the Hecke algebra
$\H(\W)$, where $\phi_{\w, \w}^d$ is identified with $T_d$.
\end{thm}

\begin{proof}
See \cite[Theorem 2.2.4]{G15} for part (i), and 
\cite[Proposition 2.2.5]{G15} for parts (ii) and (iii).
\end{proof}

Note again that parts (ii) and (iii) of Theorem \ref{Thm:1.3.6} only
apply if $n \ge r$.

\subsection{Quantum groups and tensor space}\label{ss:1.4}

The affine $q$-Schur algebras are closely related to certain quantum
groups (Hopf algebras).  The following Hopf algebra is crucial for our
purposes.

\begin{defn}\label{Def:1.4.1}
The associative, unital algebra $\ugln$ over $\qv$ is given by
generators 
$$
E_i, F_i\ \ (1 \le i \le n-1);\quad  K_i, K_i^{-1}\ \  (1 \le i \le n)
$$ 
subject to the following relations: 
\begin{align}
K_i K_j &= K_j K_i, \tag{Q1}\\
K_i K_i^{-1} &= K_i^{-1} K_i = 1, \tag{Q2}\\
K_i E_j &= v^{{\epsilon^{+}}(i, j)} E_j K_i, \tag{Q3}\\
K_i F_j &= v^{{\epsilon^{-}}(i, j)} F_j K_i, \tag{Q4}\\
E_i F_j - F_j E_i &= \delta_{ij} 
\frac{K_i K_{i + 1}^{-1} - K_i^{-1} K_{i+1}}{v - v^{-1}}, \tag{Q5}\\
E_i E_j &= E_j E_i \quad \text{ if $i$ and $j$ are not adjacent,} \tag{Q6}\\
F_i F_j &= F_j F_i \quad \text{ if $i$ and $j$ are not adjacent,}
\tag{Q7}\\
E_i^2 E_j - (v + v^{-1}) &E_i E_j E_i + E_j E_i^2 = 0 \quad 
\text{ if $i$ and $j$ are adjacent,} \tag{Q8}\\
F_j^2 F_i - (v + v^{-1}) &F_j F_i F_j + F_i F_j^2 = 0 \quad
\text{ if $i$ and $j$ are adjacent.} \tag{Q9}
\end{align}
Here, we regard $i$ and $j$ as ``adjacent'' if $i$ and $j$ index
adjacent nodes in the Dynkin diagram of type $A_{n-1}$. In the
relations, $i$ and $j$ vary over all values of the indices for which
the relation is defined. Also,
$$
\epsilon^+ (i, j) := \begin{cases}
1 & \text{ if } j = i;\\
-1 & \text{ if } \overline{j} = \overline{i-1};\\
0 & \text{ otherwise;}
\end{cases} 
$$ 
and
$$
\epsilon^- (i, j) := \begin{cases}
1 & \text{ if } \overline{j} = \overline{i-1};\\
-1 & \text{ if } j = i;\\
0 & \text{ otherwise}
\end{cases}
$$ 
where we write $\overline{a}$ for $a \in \Z$ to denote the residue
class of $a$ in the residue class ring $\Z/n\Z$. The residue class
notation has no effect in the above definition, where indices are
restricted to the range $1, \dots, n-1$. However, the notation is
important in the next two definitions.
\end{defn}

The following Hopf algebra is a quantized affine enveloping algebra 
associated with the affine Lie algebra $\widehat{\gl_n}$. 

\begin{defn}\label{Def:1.4.1A}
The associative, unital algebra $\ugn$ over $\qv$ is given by
generators $$E_i, F_i, K_i, K_i^{-1}$$ (where $1 \le i \le n$) subject
to relations (Q1) to (Q9) of Definition \ref{Def:1.4.1} (reading indices
modulo $n$).

In this definition, the notion of ``adjacent'' in relations (Q1)--(Q9)
must now be interpreted in the Dynkin diagram of type
$\widehat{A}_{n-1}$. More precisely, $i$ and $j$ are to be regarded as
``adjacent'' if $i$ and $j$ index adjacent nodes in the Dynkin diagram
of type $\widehat{A}_{n-1}$. Note that $i,j$ index adjacent nodes if
and only if $i-j \equiv \pm 1 \pmod{n}$.
\end{defn}

In \cite{G15}, a larger Hopf algebra is considered. It is an extended
version of the quantized affine algebra $\ugn$ considered in
Definition \ref{Def:1.4.1A}.

\begin{defn}\label{Def:1.4.2}
The associative, unital algebra $\eugn$ over $\qv$ is given by
generators $$E_i, F_i, K_i, K_i^{-1}, R, R^{-1}$$ (where $1 \le i \le
n$) subject to relations (Q1) to (Q9) of Definition \ref{Def:1.4.1}
(reading indices modulo $n$), together with the relations
\begin{align}
R R^{-1} &= R^{-1} R = 1, \tag{Q10}\\
R^{-1} K_{\overline{i+1}} R &= K_i, \tag{Q11}\\
R^{-1} K_{\overline{i+1}}^{-1} R &= K_i^{-1}, \tag{Q12}\\
R^{-1} E_{\overline{i+1}} R &= E_i, \tag{Q13}\\
R^{-1} F_{\overline{i+1}} R &= F_i. \tag{Q14}
\end{align}
\end{defn}

The following result was proved in \cite[Theorem 3.1.10]{G15}. 

\begin{thm}\label{Thm:1.4.3}
The algebra $\eugn$ is a Hopf algebra with multiplication $\mu$, unit
$\eta$, comultiplication $\De$, counit $\e$ and antipode $S$.
The comultiplication is defined by 
\begin{align*}
\De(1) &= 1 \otimes 1,\\
\De(E_i) &= E_i \otimes K_i K_{i+1}^{-1} + 1 \otimes E_i,\\
\De(F_i) &= K_i^{-1} K_{i+1} \otimes F_i + F_i \otimes 1,\\
\De(X) &= X \otimes X \text{\rm \ for } X \in \{K_i, K_i^{-1}, R, R^{-1}\}.
\end{align*}
The counit is defined by 
\begin{align*}
\e(E_i) &= \e(F_i) = 0,\\
\e(K_i) &= \e(K_i^{-1}) = \e(R) = \e(R^{-1}) = 1.
\end{align*}
The antipode is defined by 
\begin{align*}
S(E_i) &= -E_i K_i^{-1} K_{i+1},\\
S(F_i) &= -K_i K_{i+1}^{-1} F_i,\\
S(K_i) &= K_i^{-1},\\
S(K_i^{-1}) &= K_i,\\
S(R) &= R^{-1},\\
S(R^{-1}) &= R.
\end{align*}
The unit satisfies $\eta(1) = 1_U$.
\end{thm}
Note that the usual Hopf algebra structure on $\ugln$ and $\ugn$ is
obtained by restricting the operations of Theorem \ref{Thm:1.4.3}
above.

Let $V$ be the $\qv$-vector space with basis $\{e_t : t \in \Z\}$.
This has a natural $\eugn$-module structure as follows. 

\begin{lem}\label{Lem:1.4.4}
There is a left action of $\eugn$ on $V$ defined by the conditions 
\begin{align*}
E_i e_{t+1} &= e_t 
\text{ if } i = t \mod n, \\
E_i e_{t+1} &= 0 
\text{ if } i \ne t \mod n, \\
F_i e_{t} &= e_{t+1}
\text{ if } i = t \mod n, \\
F_i e_{t} &= 0 
\text{ if } i \ne t \mod n, \\
K_i e_t &= v e_t
\text{ if } i = t \mod n, \\
K_i e_t &= e_t
\text{ if } i \ne t \mod n, \\
R e_t &= e_{t+1}.
\end{align*}
\end{lem}

\begin{proof}
See \cite[Lemma 3.2.1]{G15}.
\end{proof}

Since $\eugn$ is a Hopf algebra, the tensor product of two
$\eugn$-modules has a natural $\eugn$-module structure via the
comultiplication $\De$.

\begin{defn}\label{Def:1.4.5}
The vector space $V^{\otimes r}$ has a natural $\eugn$-module structure
given by $u . x = \De(u)^{(r-1)} . x$.  We call this module {\em
tensor space}.  The weight $\l=(\l_1, \dots, \l_n) \in \Lambda(n, r)$
of a basis element
$$
e_{t_1} \otimes e_{t_2} \otimes \cdots \otimes e_{t_r}
$$ 
of $V^{\otimes r}$ is given by the condition 
$$
\l_i := | \{ j : t_j \equiv i \mod n \}|
$$ 
for $i = 1, \dots, n$. 
The $\l$-weight space, $V_\l$, of $V^{\otimes r}$ is the span of
all the basis vectors of weight $\l$.
\end{defn}

Henceforth, we will always assume that $q=v^2$, and regard $\qv$ as
an algebra over $\A = \Z[q,q^{-1}]$ by means of the ring homomorphism
$\A \to \qv$ such that $q \to v^2$, $q^{-1} \to v^{-2}$.

The following result about the affine $q$-Schur algebra, which will be
used frequently in the sequel, was proved in \cite[Theorem
3.4.8]{G15}.

\begin{thm}\label{Thm:1.4.6}
The quotient of $\eugn$ by the kernel of its action on tensor space is
isomorphic as a $\qv$-algebra to the algebra $\qv \otimes_\A \sqnr$.
\end{thm}

There is a corresponding result for the finite $q$-Schur algebra.
This was introduced in \cite{BLM}; see \cite{Du4} or \cite{G3} for more 
details.

\begin{thm}\label{Thm:1.4.7}
Let $V'$ be the submodule of $V$ spanned by the $e_j$ for $1 \le j \le
n$.  Then the quotient of $\ugln$ by the kernel of its action on
$(V')^{\otimes r}$ is isomorphic as a $\qv$-algebra to the algebra
$\qv \otimes_\A \bS_q(n, r)$.  We denote the corresponding epimorphism
from $\ugln$ to $\qv \otimes_\A \bS_q(n, r)$ by $\a$.
\end{thm}

\begin{defn}\label{Def:1.4.8}
For convenience of notation, we shall henceforth denote by $\svnr$ the
algebra $\qv \otimes_\A \sqnr$ and by $\bS_v(n,r)$ its finite analogue
$\qv \otimes_\A \bS_q(n,r)$. We may refer to these algebras as the 
affine $v$-Schur algebra and $v$-Schur algebra, respectively. 
\end{defn}

It will be useful in the sequel to consider the weight spaces of
$\bS_v(n, r)$ as right $\qv \otimes_\A \H(\mathcal{S}_r)$ modules.  The
following result is useful in such a context.

\begin{lem}\label{Lem:1.4.8}
Let $1 \le i_1 \le i_2 \le \cdots \le i_r \le n$, and let $\l \in
\Lambda(n, r)$ be such that $\l_j$ is the number of occurrences of $j$
in the sequence $(i_1, i_2, \ldots, i_r)$.  Then the $\l$-weight space
of $V^{\prime \otimes r}$ is generated as a right
$\qv \otimes_\A \H(\mathcal{S}_r)$-module by the element
$$
e_{i_1} \otimes e_{i_2} \otimes \cdots \otimes e_{i_r}.
$$
\end{lem}

\begin{proof}
This is a well known result, which can be seen for example by using the
definition of $\bS_v(n, r)$ together with the isomorphism, given in
\cite{Du4}, between tensor-space and Dipper and James' ``$q$-tensor
space'' (see \cite{DJ2}).
\end{proof}

Although $(V')^{\otimes r}$ is not a $\eugn$-module, we have the
following

\begin{lem}\label{Lem:1.4.9}
The action of $u \in \eugn$ on $V^{\otimes r}$ is determined by its
action on the subspace $(V')^{\otimes r}$.
\end{lem}

\begin{proof}
This is part of \cite[Proposition 3.2.5]{G15}.
\end{proof}

\subsection{Lusztig's approach}\label{ss:1.5}
In \S 1.5 we review the approach to the affine $q$-Schur algebra used
by Lusztig \cite{L11}, McGerty \cite[\S2]{Mc} and others.

Let $V_\epsilon$ be a free rank $r$ module over
$\mathbf{k}[\epsilon,\epsilon^{-1}]$, where $\mathbf{k}$ is a finite
field of $q$ elements, and $\epsilon$ is an indeterminate.

Let $\mathcal F^n$ be the space of $n$-step periodic lattices,
i.e. sequences $\mathbf L = (L_i)_{i \in \zed}$ of lattices in our
free module $V_\epsilon$ such that $L_i \subset L_{i+1}$, and
$L_{i-n}=\epsilon L_i$.  The group $G = \aut(V)$ acts on $\mathcal
F^n$ in the natural way. Let $\seq$ be the set of nonnegative integer
sequences $(a_i)_{i \in \zed}$, such that $a_i = a_{i+n}$ and
$\sum_{i=1}^n a_i = r$, and let $\mat$ be the set of $\zed \times
\zed$ matrices $A = (a_{i,j})_{i,j \in \zed}$ with nonnegative entries
such that $a_{i,j}= a_{i+n,j+n}$ and $\sum_{i \in [1,n], j \in \zed}
a_{i,j} = r$.  The orbits of $G$ on $\mathcal F^n$ are indexed by
$\seq$, where $\mathbf L$ is in the orbit $\mathcal F_{\mathbf a}$
corresponding to $\mathbf a$ if $a_i = \mathrm{dim}_{\mathbf
k}(L_i/L_{i-1})$.  The orbits of $G$ on $\mathcal F^n \times \mathcal
F^n$ are indexed by the matrices $\mat$, where a pair $(\mathbf L,
\mathbf L')$ is in the orbit $\mathcal O_A$ corresponding to $A$ if 
$$
a_{i,j}=\dim \left( \frac{L_i\cap L_j'}{(L_{i-1}\cap
L_j')+(L_i\cap L_{j-1}')} \right).
$$
For $A \in \mat$ let $r(A), c(A) \in \seq$ be given by $r(A)_i = 
\sum_{j \in \zed} a_{i,j}$ and $r(A)_j = \sum_{i \in \zed} a_{i,j}$.

Similarly let $\mathcal B^r$ be the space of complete periodic
lattices, that is, sequences of lattices $\mathbf L = (L_i)$ such that
$L_i \subset L_{i+1}$, $L_{i-r}=\epsilon L_i$, and
$\text{dim}_{\mathbf{k}}(L_i/L_{i-1})=1$ for all $i \in\zed$.  Let
$\mathbf b_0 = (\ldots,1,1,\ldots)$.  The orbits of $G$ on $\mathcal
B^r\times \mathcal B^r$ are indexed by matrices $A \in \mathfrak
S_{n,n,n}$ where the matrix $A$ must have $r(A)=c(A)=\mathbf b_0$.

Let $\mathfrak A_{r,q}$, $\mathfrak H_{r,q}$ and $\mathfrak T_{r,q}$
be the span of the characteristic functions of the $G$ orbits on
$\mathcal F^n \times \mathcal F^n$, $\mathcal B^r \times \mathcal B^r$
and $\mathcal F^n \times \mathcal B^r$ respectively. Convolution makes
$\mathfrak A_{r,q}$ and $\mathfrak H_{r,q}$ into algebras and
$\mathfrak T_r$ into a $\mathfrak A_{r,q}$--$\mathfrak H_{r,q}$
bimodule.  For $A \in \mat$ set 
$$ 
d_A = \sum_{i \ge k, j<l, 1\le i \le n}a_{ij}a_{kl}.
$$ 
Let $\{ e_A \colon A \in \mat\}$ be the basis of $\mathfrak
A_{r,q}$ given by the characteristic function of the orbit
corresponding to $A$, and let $\{[A] \colon A \in \mat\}$ be the basis
of $\mathfrak A_{r,q}$ given by $[A] = q^{-d_A/2} e_A$.  When $n = r$,
the subset of either basis spanned by all monomial matrices $A$ spans
$\mathfrak H_{r,q}$.

All of these spaces of functions are the specialization at $v = \sqrt
q$ of modules over $\A = \zed[v,v^{-1}]$, which we denote by
$\mathfrak A_r$, $\mathfrak H_r$ and $\mathfrak T_r$ respectively;
here $v$ is an indeterminate.

\begin{prop}[Varagnolo--Vasserot]\label{Prop:1.5.1} 
The $\A$-algebra $\mathfrak A_r$ is naturally isomorphic to the affine
$q$-Schur algebra $\svnr$ of Definition \ref{Def:1.3.4}.  Furthermore,
the isomorphism may be chosen to identify the basis of Definition
\ref{Def:1.3.5} with the basis $\{e_A : A \in \mat\}$.
\end{prop}

\begin{proof}
The necessary isomorphism is the map $\Phi$ given in 
\cite[Proposition 7.4 (a)]{VV}. 
\end{proof}

We will also need the {\it canonical basis}, $\{\{A\} \colon A \in
\mat\}$, for $\svnr$.  This is related to the basis $\{[A] \colon A
\in \mat\}$ in a unitriangular way: we have $$ \{ A \} = \sum_{A_1:
A_1 \le A} \Pi_{A_1, A}[A_1] ,$$ where $\le$ is a certain natural
partial order and the $\Pi_{A_1, A}$ are certain Laurent polynomials
(similar to the famous Kazhdan--Lusztig polynomials $P_{y, w}$ of
\cite{KL1}) satisfying $\Pi_{A, A} = 1$.  The reader is referred to
\cite[\S4]{L11} for full details, or to \cite[\S2.4]{G15} for a more
elementary construction.

An element $A \in \mat$ is said to be {\it aperiodic} if for any $p
\in \zed \backslash \{0\}$ there exists $k \in \zed$ such that $a_{k,
k+p} = 0$.  Let $\matap$ be the set of aperiodic elements in $\mat$.

\begin{thm}[Lusztig]\label{Thm:1.5.2} 
Under the identifications of Theorem \ref{Thm:1.4.6}, the subalgebra
$\ugn$ of $\eugn$ projects to the $\qv$-span of the elements
$$ 
\{ \{A\} : A \in \matap \}.
$$ 
\end{thm}

\begin{proof}
This is \cite[Theorem 8.2]{L11}. 
\end{proof}

\begin{rmk}\label{Rem:1.5.3}
Theorem \ref{Thm:1.5.2} is not true if we replace the canonical basis
by one of the other two bases so far discussed.

If we have $n > r$, elementary considerations show that every element
of $\mat$ is aperiodic.  This means that the subalgebra of $\svnr$
described in Theorem \ref{Thm:1.5.2} is in fact the whole of $\svnr$,
so that we may refer to the algebra of Theorem \ref{Thm:1.5.2} as
``the affine $q$-Schur algebra'' without confusion.  We will
concentrate on the case $n>r$ in this paper.
\end{rmk}

\subsection{Main results}\label{ss:1.6}

Our main aim is to prove the following

\begin{thm}\label{Thm:1.6.1}
Let $n > r$, and identify $\svnr$ with the quotient of $\ugn$
described in Theorem \ref{Thm:1.5.2} (see Remark \ref{Rem:1.5.3}).
Over $\qv$, the affine $v$-Schur algebra $\svnr$ is given by
generators $E_i, F_i, K_i, K_i^{-1}$ ($1 \le i \le n$) subject to
relations (Q1) to (Q9) of Definition \ref{Def:1.4.1A} (reading indices
modulo $n$), together with the relations
\begin{align}
K_1 K_2\cdots K_n &= v^r  \tag{Q15}\\
(K_i - 1)(K_i - v)(K_i - v^2)\cdots(K_i - v^r) &= 0. \tag{Q16}
\end{align}
\end{thm}

The corresponding result in finite type $A$ was proved by the first
author and A. Giaquinto.  We will appeal to it repeatedly in the
sequel.

\begin{thm}\label{Thm:1.6.2}
Identify $\bS_v(n,r)$ with the quotient of $\ugln$ described in Theorem
\ref{Thm:1.4.7}.  Over $\qv$, the $v$-Schur algebra $\bS_v(n,r)$ is
given by generators $E_i, F_i$ ($1\le i \le n-1$) and $K_i, K_i^{-1}$
($1 \le i \le n$) subject to relations (Q1) to (Q9) of Definition
\ref{Def:1.4.1}, together with the relations
\begin{align*}
K_1 K_2 \cdots K_n &= v^r \\
(K_i - 1)(K_i - v)(K_i - v^2)\cdots(K_i - v^r) &= 0.
\end{align*}
\end{thm}

\begin{proof}
This is \cite[Theorem 2.1]{DG}. 
\end{proof}

\begin{defn}\label{Def:1.6.3}
For now, we will denote by $T$ the $\qv$-algebra given by the
generators and relations of Theorem \ref{Thm:1.6.1}, and we will denote the
corresponding epimorphism from $\ugn$ to $T$ by $\be$.  The main aim
is thus to show that $T$ is isomorphic to $\svnr$.
\end{defn}

\begin{rmk}\label{Rem:1.6.4}
There is an obvious isomorphism between the algebra given by the
generators and relations of Theorem \ref{Thm:1.6.2} and the subalgebra
of $T$ generated by the images of the the $E_i$, $F_i$, $K_j$ and
$K_j^{-1}$, where $1 \le i < n$ and $1 \le j \le n$.  This means that
if a relation in $\svnr$ involving the $E_i$, $F_i$ and $K_j$ avoids
all occurrences of $E_a$ and $F_a$ for some $1 \le a \le n$, then by
Theorem \ref{Thm:1.6.2} and symmetry, the relation is a consequence of
relations (Q15) and (Q16).
\end{rmk}

The following result establishes a natural surjection from $T$ to
$\svnr$, and our main task in proving Theorem \ref{Thm:1.6.1} will be
to show that this map is an isomorphism, in other words, that
relations (Q15) and (Q16) are sufficient.

\begin{prop}\label{Prop:1.6.5}
Relations (Q15) and (Q16) of Theorem \ref{Thm:1.6.1} hold in $\svnr$, and
therefore $\svnr$ is a quotient of the algebra $T$.  (We denote the
corresponding epimorphism by $\g : T \ra \svnr$.)  
\end{prop}

\begin{proof}
Using the comultiplication on $\ugn$, it may be easily checked that 
$$
K_1 K_2 \cdots K_n - v^r
$$ 
and 
$$
(K_i - 1)(K_i - v)(K_i - v^2)\cdots(K_i - v^r)
$$ 
act as zero on the tensor space $(V')^{\otimes r}$ given in Theorem
\ref{Thm:1.4.7}.  The result now follows from Lemma \ref{Lem:1.4.9}.
\end{proof}

\begin{rmk}\label{Rem:1.6.6}
For later reference, we note that the maps $\alpha, \beta, \gamma$
respectively from Theorem~\ref{Thm:1.4.7}, Remark~\ref{Rem:1.6.4}, and
Proposition~\ref{Prop:1.6.5} fit together into the following
commutative diagram
\begin{equation*}
\begin{gathered}
\xymatrix{
\eugn  \ar@{->>}[drr]\\
\ugn   \ar@{->>}[r]_\beta \ar@{^{(}->}[u] & T \ar@{->>}[r]_{\gamma\quad} 
       & \svnr \\
\ugln  \ar@{->>}[r]^\alpha \ar@{^{(}->}[u] & \bS_v(n,r) \ar@{^{(}->}[u] 
}
\end{gathered}
\end{equation*}
in which all horizontal maps and the diagonal one are surjections, and
all vertical maps are injections.
\end{rmk}

\section{Proof of the main results}\label{sec:2}

Most of this section is devoted to proving Theorem
\ref{Thm:1.6.1}. The final result of this section, Theorem
\ref{Thm:2.6.1}, is an equivalent formulation of Theorem
\ref{Thm:1.6.1}, compatible with Lusztig's modified form of the
quantized enveloping algebra.

\subsection{A subalgebra of $\svnr$ isomorphic to $\H(W)$}\label{ss:2.1}
A presentation for $\svnr$ in the case $n > r$ was given in
\cite[Proposition 2.5.1]{G15}, and our main strategy for proving
Theorem \ref{Thm:1.6.1} will be to adapt this presentation.

\begin{prop}\label{Prop:2.1.1}
The algebra $\svnr$ is generated by elements 
$$ \{ \phi_{\w, \w}^d : d
\in \W \} \cup \{ \phi_{\l, \w}^1 : \l \in \Lambda(n, r) \} \cup \{
\phi_{\w, \l}^1 : \l \in \Lambda(n, r) \} .
$$ 
The elements $\phi_{\w, \w}^d$ are subject to the relations of the
affine Hecke algebra of Definition \ref{Def:1.2.1} under the
identification given by Theorem \ref{Thm:1.3.6} (iii).  The generators
are also subject to the following defining relations, where $s$
denotes a generator $s_i \in \W_\l$.
\begin{align}
\phi_{\w, \l}^1 \phi_{\mu, \w}^1 &= \d_{\l, \mu} \sum_{d \in \W_\l}
\phi_{\w, \w}^d,  \tag{Q17}\\
\phi_{\w, \w}^s \phi_{\w, \l}^1 &= q \phi_{\w, \l}^1, \tag{Q18} \\
\phi_{\l, \w}^1 \phi_{\w, \w}^s &= q \phi_{\l, \w}^1. \tag{Q19}
\end{align}
\end{prop}

A key step in understanding the structure of the algebra $T$ of
Definition \ref{Def:1.6.3} is locating within it a subalgebra
isomorphic to the affine Hecke algebra $\H(\W)$.  Theorem
\ref{Thm:1.3.6} (iii) shows that this can be done for the algebra
$\svnr$, and we now review how this works in terms of endomorphisms of
tensor space.  Recall the definition of weight space from Definition
\ref{Def:1.4.5}, and the definition of the weight $\w$ from Definition
\ref{Def:1.3.1}.

\begin{defn}\label{Def:2.1.2}
For each $1 \le i < r$, let $\t(T_{s_i}) : V_\w \rightarrow V_\w$ be
the endomorphism corresponding to the action of $v F_i E_i - 1 \in
\eugn$.  Similarly, let $\t(T_{\rho^{-1}})$ be the endomorphism
corresponding to
$$
F_n F_{n-1} \cdots F_{r+1} R,
$$ 
and let $\t(T_\rho)$ be the endomorphism corresponding to
$$
E_r E_{r+1} \cdots E_{n-1} R^{-1}.
$$
\end{defn}

\begin{lem}\label{Lem:2.1.3}
The endomorphisms $\t(T_w)$ defined above (for $w \in \{s_i: 1 \le i <
r\} \cup \{\rho, \rho^{-1}\}$) satisfy the relations of Lemma
\ref{Lem:1.2.2} (after replacing $T_w$ by $\t(T_w)$).
\end{lem}

\begin{proof}
Using the epimorphism $\a' : \ugln \twoheadrightarrow \bS_v(n, r)$
studied in \cite{BLM}, \cite{Du4}, \cite{G3}, one finds that the
action of $\t(T_{s_i})$ on $V_\w$ in the case where $i \ne r$
corresponds to the action of $\phi_{\w, \w}^{s_i} \in \bS_v(n, r)$.
(Recall from Theorem \ref{Thm:1.4.7} that $\bS_v(n, r)$ is the
quotient of $\ugln$ by the annihilator of $V_n^{\otimes r}$.)  This
proves relations $(1')$, $(2')$ and $(3')$ of Lemma \ref{Lem:1.2.2}.

The effect of $\t(T_\rho)$ on $V_\w$ is 
$$ 
\t(T_\rho)(e_{i_1} \otimes \cdots \otimes e_{i_r}) = e_{j_1}
\otimes \cdots \otimes e_{j_1} ,
$$
where $j_t = i_t - 1 \mod r$.  The effect of $\t(T_{\rho^{-1}})$ on
$V_\w$ is the inverse of this action.  The proof of relations $(4')$
and $(5')$ now follow by calculation of the action of $v F_i E_i - 1$
on $V_\w$ using the comultiplication.
\end{proof}

\begin{rmk}\label{Rem:2.1.4}
Definition \ref{Def:2.1.2} and Lemma \ref{Lem:2.1.3} are very similar
to \cite[Definition 3.3.1]{G15} and \cite[Lemma 3.3.2]{G15},
respectively.  They are included here because \cite[Definition
  3.3.1]{G15} contains an incorrect definition for $\t(T_{s_r})$.
\end{rmk}

\begin{lem}\label{Lem:2.1.5}
Define $\t(T_{s_r}) := \t(T_\rho) \t(T_{s_1}) \t(T_{\rho^{-1}})$.
Then the map taking $\t(T_w)$ to $T_w$ (where $w \in \{s_i: 1 \le i
\le r\} \cup \{\rho, \rho^{-1}\}$) extends uniquely to an isomorphism
of algebras between $\H(\W)$ and the algebra $\t(\H)$ generated by the
endomorphisms $\t(T_w)$.
\end{lem}

\begin{proof}
This follows from Lemma \ref{Lem:1.2.2} and the argument given in
\cite{G15}, namely \cite[Lemma 3.3.3, Lemma 3.3.4]{G15}.
\end{proof}

For later purposes, it will be convenient to have versions of the
above results that do not make reference to the grouplike elements $R$
and $R^{-1}$.  The following lemma is the key to the necessary
modifications.  (Recall that $n \ge r+1$ by assumption.)

\begin{lem}\label{Lem:2.1.6}
Let $e \in V_\omega$.  Then we have 
$$
R . e = (F_1 F_2 \cdots F_r) . e
$$ 
and 
$$
R^{-1} . e = (E_{r-1} E_{r-2} \cdots E_1) E_n . e.
$$
\end{lem}

\begin{proof}
It is enough to consider the case where 
$$ 
e = e_{i_1} \otimes e_{i_2} \otimes \cdots \otimes e_{i_r}
$$ 
is a basis element, and this turns out to be a straightforward 
exercise using the comultiplication in $\eugn$.
\end{proof}

\begin{prop}\label{Prop:2.1.7}
For each $1 \le i < r$, let $\t'(T_{s_i}) : V_\w \rightarrow V_\w$ be
the endomorphism corresponding to the action of $v F_i E_i - 1 \in
\ugn$.  Similarly, let $\t'(T_{\rho^{-1}}) $ be the endomorphism
corresponding to 
$$ 
(F_n F_{n-1} \cdots F_{r+1}) (F_1 F_2 \cdots F_r),
$$ 
and let $\t'(T_\rho)$ be the endomorphism corresponding to 
$$ 
(E_r E_{r+1} \cdots E_{n-1}) (E_{r-1} E_{r-2} \cdots E_1) E_n .
$$ 
Then, after replacing $T_w$ by $\t'(T_w)$, these endomorphisms
satisfy the relations of Lemma \ref{Lem:1.2.2}.
\end{prop}

\begin{proof}
Combine Lemma \ref{Lem:2.1.6} with Lemma \ref{Lem:2.1.3}.
\end{proof}

\subsection{Weight space decomposition of $T$} \label{ss:2.2}

An important property of the algebra $T$ is that it possesses a
decomposition into left and right weight spaces, similar to that
enjoyed by the ordinary and affine $q$-Schur algebras.

\begin{defn}\label{Def:2.2.1}
An element $t \in T$ is said to be of {\it left weight $\l \in
\Lambda(n, r)$} if for each $i$ with $1 \le i \le n$ we have 
$$
\be(K_i) . t = v^{\l_i} t .
$$
 where $\be$ is the map defined in Definition \ref{Def:1.6.3}.
There is an analogous definition for elements of {\it right weight
$\l$}.  The left (respectively, right) {\it $\l$-weight space} of $T$
is the $\qv$-submodule spanned by all elements of left (respectively,
right) weight $\l$.
\end{defn}

\begin{defn}\label{Def:2.2.2}
For each $\l \in \Lambda(n, r)$, define the idempotent element $1_{\l}
\in T$ by the image of $1_\l \in \bS_v(n, r)$ under the canonical
inclusion map from Remark \ref{Rem:1.6.4}. Here the $1_\l$ are the
idempotents which were defined in \cite[(3.4)]{DG}. The sum of the
$1_\l$, as $\l$ varies over $\Lambda(n,r)$, is 1. Moreover, $1_\l 1_\m
= 0$ for $\l \ne \m$, {\em i.e.}\ the idempotents are pairwise
orthogonal.
\end{defn}

\begin{prop}\label{Prop:2.2.3}
The algebra $T$ is the direct sum of its left $\l$-weight spaces, and
the nonzero $\l$-weight spaces are indexed by the elements of
$\Lambda(n, r)$.
\end{prop}

\begin{proof}
Thanks to the above orthogonal decomposition of the identity in $T$,
there is a direct sum decomposition
\[ \textstyle
T = \bigoplus_{\l\in \Lambda(n,r)} 1_\l T.
\]
Moreover, in $\bS_v(n,r)$ we have the identity
\[
   \alpha(K_i) 1_\l = \l_i 1_\l  \qquad(i = 1, \dots, n)
\]
from \cite[Proposition 8.3(a)]{DG}, where $\alpha$ is the quotient map
$\ugln \to \bS_v(n,r)$ of Theorem \ref{Thm:1.4.7}. Now it follows from
the embedding of Remark \ref{Rem:1.6.4}, or more precisely from the
commutativity of the diagram in Remark \ref{Rem:1.6.6}, that
\[
   \beta(K_i) 1_\l = \l_i 1_\l  \qquad(i = 1, \dots, n)
\]
holds in the algebra $T$. Thus it follows that $\beta(K_i) v = \l_i v$
for all $i=1,\dots, n$ and all $v \in 1_\l T$. This proves that $1_\l
T$ is the $\l$-weight space in $T$. 
\end{proof}

{\em For simplicity's sake, we will write $E_i$ in place of $\be(E_i)$
and $F_i$ in place of $\be(F_i)$ for the remainder of \S\ref{ss:2.2}.}

\begin{lem}\label{Lem:2.2.4}
{\rm(i)} In $T$ we have $K_i^{\pm 1} 1_\l = v^{\pm \l_i} 1_\l$.
\par\noindent{\rm(ii)} The idempotent $1_\l$ lies within the subalgebra of $T$
generated by the $K_i$.
\par\noindent{\rm(iii)} In $T$, the idempotent $1_\l$ coincides with the image
  of $\phi_{\l,\l}^1$ under $\beta$.
\end{lem}

\begin{proof}
Part (i) is already contained in the proof of the preceding
proposition, and part (ii) is due to the definition of $1_\l$ in
\cite{DG} as
\[\textstyle
1_\l = \left[ \begin{smallmatrix}{K_1}\\{\l_1}
\end{smallmatrix}\right] \cdots \left[ \begin{smallmatrix}{K_n}\\{\l_n}
\end{smallmatrix}\right]  .
\]
where $\left[ \begin{smallmatrix}{K_i}\\{t}
\end{smallmatrix}\right]  = \prod_{s=1}^t \frac{K_iv^{-s+1} -
  K_i^{-1}v^{s-1}}{v^s-v^{-s}}$.

Part (iii) is a consequence of the remarks preceding \cite[Lemma
2.9]{G3} combined with \cite[Corollary 2.10]{G3}.
\end{proof}

\begin{defn}\label{Def:2.2.5}
For each $i$ with $1 \le i \le n$, let $\a_i = ((\a_i)_1, \dots,
(\a_i)_n)$ be the $n$-tuple of integers given by
$$ (\a_i)_j = 
 \begin{cases} 
   1 & \text{ if } j \equiv i \mod n,\\ 
  -1 & \text{ if } j \equiv i+1 \mod n,\\ 
   0 & \text{ otherwise.} 
 \end{cases}
$$
\end{defn}

The following identities will be used frequently in the sequel, often
without explicit reference. In these identities, it will be convenient
to regard a weight $\lambda = (\lambda_1, \dots, \lambda_n)$ as an
infinite periodic sequence of integers, indexed by $\Z$, by setting
$\lambda_j$ for any $j \in \Z$ to the corresponding value $\lambda_i$
such that $1 \le i \le n$ and $j \equiv i$ mod $n$.  

\begin{lem}\label{Lem:2.2.6}  
Let $\l \in \Lambda(n, r)$, extended to an infinite periodic sequence
as above. The following identities hold in $T$:

\par\noindent{\rm (i)}
For any $1 \le i \le n$, we have $$
E_i 1_\l = \begin{cases}
1_{\lambda + \a_i} E_i & \text{ if } \lambda_{i+1} > 0;\\
0 & \text{ otherwise.}
\end{cases}
$$

\par\noindent{\rm (ii)}
For any $1 \le i \le n$, we have $$
F_i 1_\l = \begin{cases}
 1_{\lambda - \a_i} F_i & \text{ if } \lambda_{i} > 0;\\
 0 & \text{ otherwise.}
\end{cases}
$$
\end{lem}

\begin{proof}
By Remark \ref{Rem:1.6.4} and Lemma \ref{Lem:2.2.4}, it is enough to
check that both sides of each identity agree after projection to
$\svnr$.  By Theorem \ref{Thm:1.4.6}, it is enough to check that both
sides of each identity agree in their action on tensor space, which is
a routine calculation.
\end{proof}

The following lemma will be used extensively in the sequel.  We will
sometimes refer to it as the {\it cancellation principle} for $T$.

\begin{lem}\label{Lem:2.2.7} 
Let $c \ge 1$. Let $\l \in \Lambda(n, r)$, extended to an infinite
periodic sequence as above.  The following identities hold in $T$:

\par\noindent{\rm (i)}
For each $1 \le i \le n$ with $\l_i = 0$, there exists a nonzero 
element $z \in \A$ such that $$
F_i^c E_i^c 1_\l = \begin{cases}
z 1_{\lambda} & \text{ if } \lambda_{i+1} \ge c \, ;\\
0 & \text{ otherwise.}
\end{cases}
$$  Furthermore, if $c = \l_{i+1} = 1$ then $z = 1$.

\par\noindent{\rm (ii)}
For each $1 \le i \le n$ with $\l_{i+1} = 0$, there exists a nonzero 
$z' \in \A$ such that $$
E_i^c F_i^c 1_\l = \begin{cases}
z' 1_{\lambda} & \text{ if } \lambda_{i} \ge c\, ;\\
0 & \text{ otherwise.}
\end{cases}
$$  Furthermore, if $c = \l_i = 1$, then $z' = 1$.
\end{lem}

\begin{proof}
By the formulas in \cite[3.1.9]{L7} we have the following identities
in $\ugn$:
\[
\begin{aligned}
E_i^{(c)} F_i^{(c)} &= \sum_{t\ge0} F_i^{(c-t)} \prod_{s=1}^t
\frac{v^{2t-2c-s+1} \widetilde{K}_i - v^{-2t+2c+s-1}
\widetilde{K}_i^{-1}}{v^s-v^{-s}} E_i^{(c-t)}\\
F_i^{(c)} E_i^{(c)} &= \sum_{t\ge0} E_i^{(c-t)} \prod_{s=1}^t
\frac{v^{2t-2c-s+1} \widetilde{K}_i^{-1} - v^{-2t+2c+s-1}
\widetilde{K}_i}{v^s-v^{-s}} F_i^{(c-t)} 
\end{aligned} 
\]
where $\widetilde{K}_i = K_i K_{i+1}^{-1}$ and $X^{(m)}=X^m/[m]!$ for
$X=E_i,F_i$. Here $[m]$ is the quantum integer
$[m]=(v^m-v^{-m})/(v-v^{-1})$ and $[m]! = [1]\cdots[m-1][m]$ for any
$m \in \mathbb{N}$.

Since the above identities hold in $\ugn$, they hold in the quotient
$T$.  Multiply the second identity on the right by $1_\l$. By
Lemma~\ref{Lem:2.2.6} and the hypothesis $\l_i = 0$ all terms on the
right hand side will then vanish, excepting the term corresponding to
$t=c$. So we obtain the identity
\[
F_i^{(c)} E_i^{(c)} 1_\l = \prod_{s=1}^c \frac{v^{-s+1}
\widetilde{K}_i^{-1} - v^{s-1} \widetilde{K}_i}{v^s-v^{-s}} 1_\l
\]
and a similar argument with the first identity above in light of the
hypothesis $\l_{i+1} =0$ yields the identity
\[
E_i^{(c)} F_i^{(c)} 1_\l = \prod_{s=1}^c \frac{v^{-s+1}
\widetilde{K}_i - v^{s-1} \widetilde{K}_i^{-1}}{v^s-v^{-s}} 1_\l.
\]
These are identities in the quotient $T$. In fact, they hold in the
subalgebra $\bS_v(n,r)$ under the embedding of Remark~\ref{Rem:1.6.4}.
By Lemma \ref{Lem:2.2.4}(i) the above identities in $T$ take the form
\[
\begin{aligned}
F_i^{(c)} E_i^{(c)} 1_\l &= \prod_{s=1}^c \frac{v^{\l_{i+1}-\l_i-s+1}
- v^{\l_i-\l_{i+1}+s-1}}{v^s-v^{-s}} 1_\l \\ 
E_i^{(c)} F_i^{(c)} 1_\l &= \prod_{s=1}^c \frac{v^{\l_i-\l_{i+1}-s+1}
- v^{\l_{i+1}-\l_is-1} }{v^s-v^{-s}} 1_\l.
\end{aligned}
\]
Remembering that $\l_i=0$ in the first formula and $\l_{i+1}=0$ in the
second, by multiplying through by $([c]!)^2$ we obtain the desired
result, where 
\[
z = ([c]!)^2 \sqbinom{\l_{i+1}}{c}, \quad z' = ([c]!)^2 \sqbinom{\l_{i}}{c}.
\]
in terms of the standard Gaussian binomial coefficients (see {\em e.g.}
\cite[\S1.3]{L7}).  The proof is complete.
\end{proof}

\begin{defn}\label{Def:2.2.8}
Maintain the notation of Lemma \ref{Lem:2.2.7}.  Let $M$ be a monomial
in the various elements $E_i$, $F_i$ and $1_\l$ of $T$.  We call a
monomial $M'$ a {\it reduction} of $M$ if it (a) represents the same
element of $T$ as $M$ and (b) $M'$ can be obtained from $M$ by
omitting zero or more generators of $M$ of the form $1_\mu$.

A {\it distinguished term} in the algebra $T$ is an element of $T$ of
one of the following two forms:
\par\noindent{(i)}
{$E_i^c 1_\l$, where $c \ge 0$ and $\l_i = 0$;}
\par\noindent{(ii)}
{$F_i^c 1_\l$, where $c \ge 0$ and $\l_{i+1} = 0$.}

A {\it strictly distinguished monomial} in the algebra $T$ is a
monomial in the elements $F_i$, $E_i$ and $1_\l$ that can be parsed as
a word in the distinguished terms.  A reduction of a strictly
distinguished monomial is called a {\it distinguished monomial}.
\end{defn}

\begin{exam}\label{Ex:2.2.9}
The idempotents $1_\l$ are both distinguished terms and distinguished
monomials in $T$: here, we take $c = 0$.

If $\l_i=0$ and $\l_{i+1} = c$ then the element $M' = F_i^c E_i^c
1_\l$ of Lemma \ref{Lem:2.2.7} (i) is a distinguished monomial.
Indeed, it can be seen by repeated applications of Lemma
\ref{Lem:2.2.6} (i) and the fact that $1_\l$ is idempotent, that $M'$
is a reduction of the strictly distinguished monomial $M = (F_i^c
1_{\l + c\a_i}) (E_i^c 1_\l)$.  (To verify that $M$ is strictly
distinguished, one must note that $(\l + c\a_i)_{i+1} = 0$.)
Furthermore, by Lemma \ref{Lem:2.2.7} and the fact that $\l_{i+1} \ge
c$, we see that $M = M'$ is a nonzero element of $T$.

Similarly, if $\l_i = c$ and $\l_{i+1} = 0$ then the element $E_i^c
F_i^c 1_\l$ of Lemma \ref{Lem:2.2.7} (ii) is a distinguished monomial.
\end{exam}

In the sequel, we will make use of various automorphisms of $T$; part
(i) below may be used without explicit comment.

\begin{prop} \label{Prop:2.2.10}
{\rm (i)}
There is a unique automorphism $\nu$ of $T$ of order $n$ satisfying 
\begin{align*}
\nu(E_i) &= E_{i+1},\\
\nu(F_i) &= F_{i+1} \text{ and} \\
\nu(K_i^{\pm 1}) &= K_{i+1}^{\pm 1}, 
\end{align*}
for all $1 \le i \le n$, and reading subscripts modulo $n$.  Let $\l
\in \Lambda(n, r)$ and define
$$
\l_+ = (\l_n, \l_1, \l_2, \l_3, \ldots, \l_{n-1}).
$$  
Then $\nu(1_\l) = 1_{\l_+}$.

\par\noindent{\rm (ii)}
There is a unique anti-automorphism $\s$ of $T$ satisfying 
\begin{align*}
\s(E_i) &= F_i, \\
\s(F_i) &= E_i \text{ and} \\
\s(K_i^{\pm 1}) &= K_i^{\pm 1}, 
\end{align*}
for all $1 \le i \le n$.  The anti-automorphism $\s$ fixes all
elements $1_\l \in T$.
\end{prop}

\begin{proof}
For (i), we note that there is an automorphism of $\ugn$ corresponding
to $\nu$, that in addition fixes the elements $R^{\pm 1}$; this can be
verified by checking the defining relations for $\ugn$.  Since this
automorphism preserves setwise the set of relations (Q15) and (Q16) in
$T$, we obtain an automorphism of $T$ as claimed; it is unique because
we have given its effect on a generating set (see Theorem
\ref{Thm:1.6.1}).  The last claim of (i) follows from the relationship
between the $K_i$ and $1_\l$; see for example \cite[Corollary
2.10]{G3}.

The same line of argument can be used to prove (ii).
\end{proof}

\begin{lem}\label{Lem:2.2.11}
Let $M = t_1 t_2 \cdots t_k$ be a strictly distinguished monomial with
distinguished terms $t_i$.  Then $M \ne 0$ if and only if the
following two conditions hold:
\par\noindent{\rm (i)} each term $t_i$ is nonzero;
\par\noindent{\rm (ii)} for each $1 \le i < k$, there exists $\l \in \Lambda(n, r)$ 
such that $t_i = t_i 1_\l$ and $t_{i+1} = 1_\l t_{i+1}$.
\end{lem}

\begin{proof}
Condition (i) is clearly necessary for $M$ to be nonzero.  To see the
necessity of condition (ii), recall from Lemma \ref{Lem:2.2.6} that
for each term $t_i$, there exist $\l, \mu \in \Lambda(n, r)$ such that
$t_i = 1_\l t_i = t_i 1_\mu$.

We now check sufficiency.  It will be enough to show that $\s(M) M \ne
0$, where $\s$ is as in Proposition \ref{Prop:2.2.10}.  This follows
from Lemma 2.2.7.  Indeed, the hypotheses $\l_{i+1} \ge c$ or $\l_i
\ge c$ follow from condition (i) above, and condition (ii) above
implies that if $t_i 1_\l = t_i$ then we have
\begin{align*}
\s(t_{i+1}) \s(t_i) t_i t_{i+1} &= z'' \s(t_{i+1}) 1_\l t_{i+1} \\ 
        &= z'' \s(t_{i+1}) t_{i+1}, 
\end{align*}
where $z''$ is equal either to $z$ or to $z'$ as in Lemma
\ref{Lem:2.2.7}.  There is a unique $\mu \in \Lambda(n, r)$ such that
$M = M 1_\mu$, and an induction now shows that $\s(M) M$ is a nonzero
scalar multiple of $1_\mu$, completing the proof.
\end{proof}

\subsection{A subalgebra of $T$ isomorphic to $\svnr$}\label{ss:2.3}

The aim of \S\ref{ss:2.3} is to show that the relations satisfied by
the endomorphisms of Proposition \ref{Prop:2.1.7} are in fact
consequences of the defining relations (Q15) and (Q16) of the algebra
$T$.  In this section, we may abuse notation by identifying elements
$u$ of $\ugn$ with their images $\be(u)$ in $T$ (see Definition
\ref{Def:1.6.3}).

Recall from Remark \ref{Rem:1.6.4} that there is a natural subalgebra
of $T$ isomorphic to the ordinary $v$-Schur algebra, $\bS_v(n, r)$.
Using this fact, we can make the following

\begin{defn}\label{Def:2.3.1}
For each $1 \le i < r$, define elements of $T$ by
\begin{align*}
\z(T_{s_i}) &= (v F_i E_i - 1) 1_\w,\\
\z(T_{\rho^{-1}}) &= \left(
(F_n F_{n-1} \cdots F_{r+1}) (F_1 F_2 \cdots F_r) \right) 1_\w,\\
\z(T_\rho) &= \left(
(E_r E_{r+1} \cdots E_{n-1}) (E_{r-1} E_{r-2} 
\cdots E_1) E_n \right) 1_\w .
\end{align*}
\end{defn}

\begin{rmk}\label{Rem:2.3.2}
It follows from repeated applications of Lemma \ref{Lem:2.2.6} that
each element $\z(T_w)$ given in Definition \ref{Def:2.3.1} has the
property that $1_\w \z(T_w) = \z(T_w)$.
\end{rmk}

\begin{lem}\label{Lem:2.3.3}
The expressions given for $\z(T_{\rho^{-1}})$ and $\z(T_{\rho})$ are
distinguished monomials.
\end{lem}

\begin{proof}
This is a routine exercise, in which the hypothesis that $n > r$ plays
an important part.
\end{proof}

\begin{lem}\label{Lem:2.3.4}
The following identities hold in $T$: 

\par\noindent{\rm(i)} $\z(T_{\rho^{-1}}) = \left(
F_n (F_1 F_2 \cdots F_{r-2} F_{r-1}) (F_{n-1} F_{n-2} \cdots F_{r+1} F_r)
\right) 1_\w$;

\par\noindent{\rm(ii)} $\z(T_\rho) = \left(
(E_r E_{r-1} \cdots  E_2 E_1) (E_{r+1} E_{r+2} \cdots E_{n-1} E_n)
\right) 1_\w$.
\end{lem}

\begin{proof}
Equation (i) (respectively, (ii)) follows by applying repeated
commutations between the generators $F_i$ (respectively, $E_i$).
\end{proof}

\begin{lem}\label{Lem:2.3.5}
The following identities hold in $T$:
\par\noindent{\rm (i)} $\z(T_{\rho^{-1}}) \z(T_{\rho}) = 1_\w;$
\par\noindent{\rm (ii)} $\z(T_{\rho}) \z(T_{\rho^{-1}}) = 1_\w.$
\end{lem}

\begin{proof}
Let $\w' \in \Lambda(n, r)$ be the weight $\w' = (0, 1, \ldots, 1, 0,
\ldots, 0),$ where the occurrences of $1$ appear in positions $2, 3, 4,
\ldots, r+1 \le n$.  Then it follows from Lemma \ref{Lem:2.3.4} that
$$
\z(T_\rho) = 1_\w \left(
(E_r E_{r-1} \cdots  E_2 E_1) 1_{\w'} (E_{r+1} E_{r+2} \cdots E_{n-1} E_n)
\right) 1_\w
$$ 
and it follows from Definition \ref{Def:2.3.1} that 
$$
\z(T_{\rho^{-1}}) = 1_\w \left(
(F_n F_{n-1} \cdots F_{r+1}) 1_{\w'} (F_1 F_2 \cdots F_r) \right) 1_\w.
$$
We will prove (i), and (ii) follows by a similar argument.

To prove (i), we first show that 
$$
1_{\w'} (F_1 F_2 \cdots F_r) 1_\w (E_r E_{r-1} \cdots  E_2 E_1) 1_{\w'}
= 1_{w'}.
$$ 
The left hand side of the equation is readily checked to be a good
monomial, and then the equation follows by repeated applications of
the $c = 1$ case of the cancellation principle (Lemma
\ref{Lem:2.2.7}), starting in the middle of the equation (i.e., with
$F_r 1_\w E_r$).  A similar argument shows that
$$ 
1_\w (F_n F_{n-1} \cdots F_{r+1}) 1_{\w'} 
(E_{r+1} E_{r+2} \cdots E_{n-1} E_n) 1_\w = 1_\w .
$$ 
Part (i) follows by combining these last two identities.
\end{proof}

\begin{lem}\label{Lem:2.3.6}
Let $1 < i < r$, and let 
$$
M = (E_{r-1} E_{r-2} \cdots E_1) (E_{r+1} E_{r+2} \cdots E_n) 1_\w.
$$ 
Then the identity $(v F_{i-1} E_{i-1} - 1) M = M (v F_i E_i - 1)$
holds in $T$.  
\end{lem}

{\em Note.} Notice that both sides of the identity have right weight
$\w$.

\begin{proof}
Since the identity involves no occurrences of $E_r$ or $F_r$, Remark
\ref{Rem:1.6.4} applies.  More precisely, after applying a suitable
symmetry of the Dynkin diagram, we see that it suffices to prove the
identity
$$
(v F_{n-r+i-1} E_{n-r+i-1} - 1)M' = M' (v F_{n-r+i} E_{n-r+i} - 1)
$$ 
in the ordinary $v$-Schur algebra, where indices are read modulo $n$,
and we have
$$
M' = (E_{n-1} E_{n-2} \cdots E_{n - r + 1})(E_1 E_2 \cdots E_{n-r}) 1_\w,
$$ 
and 
$$
\w' = (0, 0, \ldots, 0, 1, 1, \ldots, 1),
$$ 
where $1$ occurs $r$ times in $\w'$.  

By Theorem \ref{Thm:1.4.7}, it suffices to show that both sides of the
identity act in the same way on tensor space $V^{\prime \otimes r}$,
and because both sides of the identity have right weight $\w'$, it is
enough to check this on the $\w'$-weight space.  By Lemma
\ref{Lem:1.4.8}, it is enough to check that each side of the identity
acts the same on the element
$$ 
e_{\w'} = e_{n-r+1} \otimes e_{n-r+2} \otimes \cdots \otimes e_n .
$$

Fix $j$ with $1 \le j < r$, and let $e_{j, \w'}$ be the tensor
obtained by exchanging the occurrences of $e_{n-r+j}$ and
$e_{n-r+j+1}$ in $e_{\w'}$.  Using the comultiplication, it is a
routine calculation to show that 
$$ 
(F_{n-r+j} E_{n-r+j}) e_{\w} = e_{j, \w'} + v^{-1} e_{\w'} ,
$$ 
and it is immediate from this that 
$$
(v F_{n-r+j} E_{n-r+j} - 1) e_\w = v e_{j, \w'} .
$$

Another calculation with the comultiplication shows that 
$$
(E_{n-1} E_{n-2} \cdots E_{n - r + 1}) 
(E_1 E_2 \cdots E_{n-r}) 1_{\w'} e_{\w'}
= e'_{\w'},
$$ where 
$$
e'_{\w'} = e_1 \otimes (e_{n-r+1} \otimes e_{n-r+2} \otimes \cdots \otimes
e_{n-1}).
$$  
Let $j$ be such that $1 \le j < r$.  
Acting $(v F_{n-r+j} E_{n-r+j} - 1)$ on the left, we deduce that 
$$
(v  F_{n-r+j} E_{n-r+j} - 1)
(E_{n-1} E_{n-2} \cdots E_{n - r + 1}) 
(E_1 E_2 \cdots E_{n-r}) 1_{\w'} e'_{\w'}
= v e'_{j, \w'},
$$ 
where $e'_{j, \w'}$ is obtained from $e_{j, \w}$ by exchanging the
occurrences of $e_{n-r+j}$ and $e_{n-r+j+1}$.

The result now follows after we observe that 
$$
(E_{n-1} E_{n-2} \cdots E_{n - r + 1}) 
(E_1 E_2 \cdots E_{n-r}) 1_{\w'} e_{j+1, \w'}
= 
e'_{j, \w'}.
$$
\end{proof}

\begin{cor}\label{Cor:2.3.7}
If $i$ is such that $1 < i < r$, then the relation 
$$
\z(T_{s_{i-1}}) \z(T_\rho) = \z(T_\rho) \z(T_{s_i})
$$ 
holds in $T$.
\end{cor}

\begin{proof}
Observe that $E_r$ commutes with $F_i$, $F_{i-1}$, $E_i$ and
$E_{i-1}$.  The assertion now follows by left-multiplying the identity
of Lemma \ref{Lem:2.3.6} by $E_r$.
\end{proof}

The techniques of proof of Lemma \ref{Lem:2.3.6} play an important
part in the next brace of results.

\begin{lem}\label{Lem:2.3.8}
The following identities hold in $T$, where $1 < i < r$:
\par\noindent{\rm(i)} $(F_i E_i - v^{-1}) 1_\w = (E_i F_i - v)1_\w$;
\par\noindent{\rm(ii)}
$1_\w (F_n F_{n-1} \cdots F_r E_r E_{r+1} \cdots E_n - v^{-1}) = 
1_\w (E_r E_{r+1} \cdots E_n F_n F_{n-1} \cdots F_r - v)$;
\par\noindent{\rm(iii)} $(E_n F_n - v) E_1 E_n 1_\w = E_1 E_n (E_1 F_1  - v) 1_\w;$
\par\noindent{\rm(iv)} $1_\w (F_{r-1} E_{r-1} - v^{-1}) E_r E_{r-1} \cdots E_2 E_1
= 1_\w E_r E_{r-1} \cdots E_2 E_1 (F_r E_r - v^{-1}).$
\end{lem}

{\em Note.}  The expressions appearing in (i) and (ii) above have both
left and right weight equal to $\w$.

\begin{proof}
We omit the proof of (i), because it is similar to, but easier than, the
proof of (ii).

To prove (ii), it is enough, by symmetry of the defining relations of $T$,
to prove the identity 
$$
1_{\w^-} (F_{n-1} \cdots F_{r-1} E_{r-1} \cdots E_{n-1} - v^{-1}) = 
1_{\w^-} (E_{r-1} \cdots E_{n-1} F_{n-1} \cdots F_{r-1} - v),
$$ where 
$$
\w^- = (\underbrace{1, 1, \ldots, 1}_{r-1}, 
\underbrace{0, 0, \ldots, 0}_{n-r}, 1).
$$ 
This can be regarded as an identity in $\bS_v(n, r)$.  By Lemma
\ref{Lem:1.4.8}, it is enough to check that each side of the identity
acts in the same way on the element
$$ 
e_{\w^-} = e_1 \otimes e_2 \otimes \cdots
\otimes e_{r-2} \otimes e_{r-1} \otimes e_n .
$$ 
A calculation shows that each side of the identity acts on $e_{\w^-}$ 
to give 
$$ 
e_1 \otimes e_2 \otimes \cdots \otimes e_{r-2} \otimes e_n \otimes e_{r-1}.
$$

To prove (iii), it is enough by symmetry of the defining relations to prove 
$$
(E_1 F_1  - v) E_2 E_1 1_{\w^+} = E_2 E_1 (E_2 F_2 - v) 1_{\w^+},
$$ 
where 
$$
\w^+ = (0, \underbrace{1, 1, \ldots, 1}_r, 
\underbrace{0, 0, \ldots, 0}_{n-r-1}).
$$ 
By Lemma \ref{Lem:1.4.8}, it is enough to show that both sides of
the identity act in the same way on
$$
e_{\w^+} = e_2 \otimes e_3 \otimes \cdots \otimes e_{r+1}.
$$  
A calculation shows that each side sends $e_{\w^+}$ to 
$$
e_2 \otimes e_1 \otimes e_4 \otimes e_5 \otimes \cdots \otimes e_{r+1}.
$$

For (iv), observe that both sides of the identity have right weight
$\w^+$, as defined above.  Since (iv) can be regarded as an identity
in $\bS_v(n, r)$, Lemma \ref{Lem:1.4.8} applies and it is enough to
check that both sides of the identity have the same effect on
$e_{\w^+}$.  A calculation shows that both sides of the identity send
$e_{\w^+}$ to
$$ 
e_1 \otimes e_2 \otimes \cdots \otimes e_{r-3} \otimes e_{r-2}
\otimes e_r \otimes e_{r-1} .
$$
\end{proof}

\begin{defn}\label{Def:2.3.9}
We define $\z(T_{s_r})$ to be the element of $T$ given by 
$$
1_\w (F_n F_{n-1} \cdots F_r E_r E_{r+1} \cdots E_n - v^{-1}).
$$
\end{defn}

\begin{lem}\label{Lem:2.3.10}
The following identities hold in $T$: 
\par\noindent{\rm (i)} $\z(T_\rho) \z(T_{s_r}) = \z(T_{s_{r-1}}) \z(T_\rho);$
\par\noindent{\rm (ii)} $\z(T_\rho) \z(T_{s_1}) = \z(T_{s_r}) \z(T_\rho).$
\end{lem}

\begin{proof}
We prove (i) first.  Using Lemma~\ref{Lem:2.3.4}(ii), it is enough to
prove that the expressions
$$
M_1 = 1_\w (E_r E_{r-1} \cdots E_1)
(E_{r+1} E_{r+2} \cdots E_n) (F_n F_{n-1} \cdots F_{r+1} F_r) (E_r E_{r+1}
\cdots E_n)
$$ 
and 
$$
M_2 = 1_\w (F_{r-1} E_{r-1}) (E_r E_{r-1} \cdots E_1) (E_{r+1} E_{r+2} 
\cdots E_n)
$$ 
are equal.

Using Lemma \ref{Lem:2.2.6} repeatedly, and the notation of the proof
of Lemma \ref{Lem:2.3.6}, we find that
$$
M_1 = 1_\w E_r E_{r-1} \cdots E_1
(1_{\w^+}
E_{r+1} E_{r+2} \cdots E_n F_n F_{n-1} \cdots F_{r+1}
1_{\w^+})
F_r E_r E_{r+1}
\cdots E_n
,
$$ 
and repeated applications of the cancellation principle (Lemma
\ref{Lem:2.2.7}) show that the parenthetic expression shown is simply
equal to $1_{\w^+}$.  By Lemma \ref{Lem:2.2.6} again, we have
\begin{align*}
M_1 &= 1_\w E_r
E_{r-1} \cdots E_1 1_{\w^+} F_r E_r E_{r+1} \cdots E_n \\ 
&= 1_\w (E_r E_{r-1} \cdots E_1) F_r E_r (E_{r+1} E_{r+2} \cdots E_n) .
\end{align*}
Applying Lemma \ref{Lem:2.3.8} (iv) gives $$ M_1 = 1_\w F_{r-1}
E_{r-1} (E_r E_{r-1} \cdots E_1) (E_{r+1} E_{r+2} \cdots E_n) ,$$
which is $M_2$, as required.

We now turn to (ii).  By Lemma \ref{Lem:2.3.8}, parts (i) and (ii), it
is enough to show that the monomials
$$
M_3 = 1_\w 
(E_r E_{r+1} \cdots E_n)
(F_n F_{n-1} \cdots F_r)
(E_r E_{r+1} \cdots E_{n-1})
(E_{r-1} E_{r-2} \cdots E_1) E_n
$$ 
and 
$$
M_4 = 1_\w
(E_r E_{r+1} \cdots E_{n-1})
(E_{r-1} E_{r-2} \cdots E_1) E_n
(E_1 F_1)
$$ 
are equal.  Using Lemma \ref{Lem:2.2.6} and the notation of the
proof of Lemma \ref{Lem:2.3.6} again, we find that $M_3$ is equal to
$$
1_\w 
E_r E_{r+1} \cdots E_n
F_n 
(1_{\w^-} F_{n-1} F_{n-2} \cdots F_r
E_r E_{r+1} \cdots E_{n-1} 1_{\w^-})
E_{r-1} E_{r-2} \cdots E_1 E_n
.$$  
By the cancellation principle, this simplifies to 
\begin{align*}
M_3 &= 1_\w 
E_r E_{r+1} \cdots E_n
F_n 
1_{\w^-}
E_{r-1} E_{r-2} \cdots E_1 E_n \\
&= 1_\w
(E_r E_{r+1} \cdots E_{n-1})
(E_n F_n)
(E_{r-1} E_{r-2} \cdots E_1) E_n.
\end{align*}
Applying commutations yields 
$$
M_3 = 1_\w
(E_r E_{r+1} \cdots E_{n-1})
(E_{r-1} E_{r-2} \cdots E_2)
(E_n F_n E_1 E_n).
$$  
Using Lemma \ref{Lem:2.3.8} (iii), we have 
$$
M_3 = 1_\w
(E_r E_{r+1} \cdots E_{n-1})
(E_{r-1} E_{r-2} \cdots E_2)
(E_1 E_n E_1 F_1),
$$ 
which is equal to $M_4$, as desired.
\end{proof}

\begin{prop}\label{Prop:2.3.11} 
Let $\g : T \ra \svnr$ be the epimorphism of Proposition
\ref{Prop:1.6.5}.  Then $\g$ admits a right inverse: there is an
injective homomorphism
$$\iota : \svnr \ra T$$ such that $\g \circ \iota$ is the identity
homomorphism on $\svnr$.  
\end{prop}

\begin{proof}
We start by specifying $\iota$ on the subalgebra of $\svnr$ spanned by
the elements $\phi_{\w, \w}^d$, as in \S\ref{ss:2.1}.  We define
$$
\iota(\tau(T_{s_i})) := \zeta(T_{s_i})
$$ 
for $1 \le i < r$, and define 
$$
\iota(\tau(T_{\rho^{\pm 1}})) := \zeta(T_{\rho^{\pm 1}}).
$$ In these cases, $\g \circ \iota$ is the identity map by Proposition
2.1.7 and the definition of $\zeta$, so it is enough to check that the
relations of Lemma \ref{Lem:1.2.2} are satisfied in the image of
$\iota$.  The difficult cases, ($4'$) and ($5'$), follow from Lemma
\ref{Lem:2.3.5}, Corollary \ref{Cor:2.3.7} and Lemma \ref{Lem:2.3.10}.
Since those cases hold, it is enough to check cases ($1'$), ($2'$) and
($3'$) assuming that neither $s$ nor $t$ is equal to $r$; this follows
from Lemma \ref{Lem:2.1.3} and Remark \ref{Rem:1.6.4}.

It remains to check relations (Q17), (Q18) and (Q19) of Proposition
2.1.1.  Since, by Remark \ref{Rem:1.6.4}, there is a canonically
embedded copy of $\bS_v(n, r)$ in the algebra $T$ (namely, the
subalgebra generated by all $E_i$, $F_i$, $K_j^{\pm 1}$ with $i \ne
n$), we may send the elements $\phi_{\w, \l}^1$, $\phi_{\mu, \w}^1$
and $\phi_{\w, \w}^d$ (where $d$ lies in the finite symmetric group)
to the corresponding elements of $T$.  (Observe that this construction
is compatible with the definitions in the previous paragraph.)  More
explicitly, if $M$ is a polynomial in the generators $E_i$, $F_i$,
$K_j$, $K_j^{-1}$ and there are no occurrences of $E_n$ or $F_n$, then
$\a(M) \in \bS_v(n, r)$ and $\iota(\a(M)) = \be(M)$ by construction.  It
then follows that
$$
(\g \circ \iota)(\a(M)) = \g(\be(M)) = \a(M) ,
$$ as required.  Relations (Q17), (Q18) and (Q19) can be seen to hold in
$T$ (and thus in the image of $\iota$) by Theorem \ref{Thm:1.6.2} and
Remark \ref{Rem:1.6.4}.
\end{proof}

\subsection{Surjectivity of $\iota$}\label{ss:2.4}

So far we have shown that there is a monomorphism $\iota : \svnr \ra
T$. From Definition \ref{Def:1.6.3} we have a surjective map $\be :
\ugn \ra T$.  We aim in \S\ref{ss:2.4} to show that the image of $\be$
is contained in the image of $\iota$, which will complete the proof of
Theorem \ref{Thm:1.6.1}.

\begin{lem}\label{Lem:2.4.1}
The algebra $T$ is generated by the all elements of the form 
$E_i 1_\l$, $1_\l F_i$ and $1_\l$, for $1 \le i \le n$ and $\l \in
\Lambda(n, r)$.
\end{lem}

\begin{proof}
The images of the $K_i$ and $K_i^{-1}$ in $T$ are linear combinations
of the elements $1_\l$, as they are in the ordinary $v$-Schur algebra
(see \cite[Corollary 2.10]{G3}).  The image of $E_i$ in $T$ is a
linear combination of elements $E_i 1_\l$, because the $1_\l$ form an
orthogonal decomposition of the identity (see Lemma \ref{Lem:2.2.4}).
Similarly, the image of $F_i$ in $T$ is a linear combination of
elements $1_\l F_i$.  This shows that the usual algebra generators of
$T$ lie in the span of the elements listed in the statement.
Conversely, it follows from the definitions (see Definition
\ref{Def:2.2.2}) that the elements listed lie in $T$, completing the
proof.
\end{proof}

\begin{lem}\label{Lem:2.4.2}
If $i \ne n$ and $\l \in \Lambda(n, r)$, the elements $E_i 1_\l$,
$1_\l F_i$ and $1_\l$ lie in the image of $\iota$.  
\end{lem}

\begin{proof}
The elements listed in the statement lie in the canonically embedded
copy of $\bS_v(n, r)$ in $T$.  By the construction of $\iota$ (see the
proof of Proposition \ref{Prop:2.3.11}), such elements lie in the
image of $\iota$.
\end{proof}

\begin{lem}\label{Lem:2.4.3}
{\rm (i)} The element $E_n 1_\w$ of $T$ lies in $\iota(\svnr)$.
\par\noindent{\rm (ii)} The element $1_\w F_n$ of $T$ lies in $\iota(\svnr)$.
\end{lem}

\begin{proof}
By construction of $\iota$, the element $\z(T_\rho)$ lies in the image
of $\iota$.  Since $T$ contains a canonical copy of $\bS_v(n, r)$ (see
the proof of Proposition \ref{Prop:2.3.11}), the element 
$$ (F_1 F_2 \cdots F_{r-1})(F_{n-1} F_{n-2} \cdots F_r) 1_\w
$$ 
lies in the image of $\iota$.  By multiplying these two elements we see
that 
$$
(F_1 F_2 \cdots F_{r-1})(F_{n-1} F_{n-2} \cdots F_r)
(E_r E_{r+1} \cdots E_{n-1}) (E_{r-1} E_{r-2} \cdots E_1) E_n 1_\w 
$$ 
lies in the image of $\iota$.  Applying the cancellation principle,
this latter expression simplifies to $E_n 1_\w$, completing the proof
of (i).

A similar argument using $\z(T_{\rho^{-1}})$ in place of $\z(T_\rho)$
can be used to prove (ii).
\end{proof}

Our main effort will be directed towards proving that the elements
$E_n 1_\l$ lie in the image of $\iota$.  More precisely, we will prove
that $E_n 1_\l$ lies in the ideal of $\Im(\iota)$ generated by $E_n
1_\w$.  Our argument will rely on the following technical lemma whose
proof will be deferred to \S\ref{ss:2.5}.

\begin{lem}\label{Lem:2.4.4}
Fix $\l \in \Lambda(n, r)$ such that $\l_1 > 0$.  There exists a
distinguished monomial $M$ in the generators $E_i$, $F_i$ and $1_\mu$
satisfying the following conditions:
\par\noindent{\rm (i)} $M = 1_\w M 1_\l \ne 0$;
\par\noindent{\rm (ii)} $M$ contains no occurrences of $E_n$, $F_n$, $E_1$ or 
$F_{n-1}$;
\par\noindent{\rm (iii)} all the occurrences of $F_1$ occur consecutively, as
do all the occurrences of $E_{n-1}$;
\par\noindent{\rm (iv)} there are at most $\l_1 - 1$ occurrences of $F_1$.
\end{lem}

\begin{lem}\label{Lem:2.4.5}
Let $\s$ be the antiautomorphism of $T$ given in Proposition
\ref{Prop:2.2.10} (ii), and let $M$ and $\l$ be as in Lemma
\ref{Lem:2.4.4}.  Then there is a nonzero scalar $z \in \qv$ such that
$$
E_n 1_\l = z \sigma(M) (E_n 1_\w) M.
$$ In particular, $E_n 1_\l$ lies in the ideal of $\Im(\iota)$
generated by $E_n 1_\w$.
\end{lem}

\begin{proof}
The monomial $M$ is equal to a strictly distinguished monomial 
$$
M' = t'_m t'_{m-1} \cdots t'_1 .
$$ 
After moving unnecessary idempotents in $M$ to the right using
Lemma \ref{Lem:2.2.7}, and omitting the corresponding idempotents from
the terms $t'_i$, we may assume that $M$ is of the form
$$ 
M = t_m t_{m-1} \cdots t_1 1_\l .
$$ 
We will prove by induction on $k \le m$ that 
$$
(\s(t_1) \s(t_2) \cdots \s(t_k)) E_n (t_k t_{k-1} \cdots t_1) 1_\l
$$ 
is a nonzero multiple of $E_n 1_\l$; the case $k = m$ is the
assertion of the Lemma.  The base case, $k = 0$, is trivial.

There are two cases to consider for the inductive step.  The first
case, which is easier to deal with, is that $t_k$ is of the form
$E_i^c$ for some $c > 0$.  In this case, we have
$$
\s(t_k) E_n t_k = F_i^c E_n E_i^c,
$$ 
which can be rewritten as 
$$
E_n F_i^c E_i^c 
$$ 
using the relations of $\ugn$.  (Note that we do not have $i = n$,
because $M$ does not contain occurrences of $E_n$.)  We now have
$$
F_i^c E_n t_k t_{k-1} \cdots t_1 1_\l
=
E_n F_i^c t_k 1_\mu t_{k-1} t_{k-2} \cdots t_1 1_\l
$$ 
for a suitable $\mu \in \Lambda(n, r)$.  The hypothesis $M \ne 0$
means that we have $\mu_{i+1} \ge c$, so we may apply Lemma
\ref{Lem:2.2.7} (i) to replace $F_i^c t_k 1_\mu$ by $z 1_\mu$ with $z$
nonzero.  The proof is now completed in this case by the inductive
hypothesis.

The second case is that $t_k$ is of the form $F_i^c$ for some $c > 0$.
In this case, we cannot have $i = n$ or $i = {n-1}$ because of
condition (ii) of Lemma \ref{Lem:2.4.4}.  Suppose for the moment that
$i \ne 1$.  Then the relations in $\ugn$ show that
$$ 
E_i^c E_n F_i^c = E_n E_i^c F_i^c .
$$ 
We can then proceed as in the first case to show that
\begin{align*}
E_i^c E_n t_k t_{k-1} \cdots t_1 1_\l &= E_n E_i^c t_k
1_\mu t_{k-1} t_{k-2} \cdots t_1 1_\l \\ 
&= E_n z' 1_\mu t_{k-1} t_{k-2} \cdots t_1 1_\l 
\end{align*}
where $z'$ is nonzero.  Here we have used Lemma \ref{Lem:2.2.7} (ii),
which is applicable because $M$ is nonzero and $\mu_i \ge c$.  Again,
we are done by induction in this case.

The remaining case is the possibility that $t_k = F_1^c$ for some $c >
0$.  The relations in $\ugn$ show that
$$
E_1^c E_n F_1^c = E_1^c F_1^c E_n,
$$ 
and so we have 
$$
E_1^c E_n t_k t_{k-1} \cdots t_1 1_\l
= E_1^c t_k E_n 1_\mu t_{k-1} t_{k-2} \cdots t_1 1_\l.
$$  
Because $F_1^c$ arises from a distinguished term, we have $\mu_2 = 0$.
By Lemma \ref{Lem:2.2.6} (i), we have 
$$ 
E_1^c F_1^c E_n 1_\mu = E_1^c F_1^c 1_{\mu'} E_n ,
$$ 
where again $\mu'_2 = 0$ (recall that $n \ge 3$).  Hence 
$$ 
E_1^c E_n t_k t_{k-1} \cdots t_1 1_\l = z' 1_{\mu'} E_n t_{k-1}
t_{k-2} \cdots t_1 1_\l .
$$ Furthermore, $z'$ is nonzero.  To see why, we recall that by
condition (ii) of Lemma \ref{Lem:2.4.4}, there are no occurrences of
$E_1$ or $E_n$ or $F_n$ in $M$ and that by condition (iii), all the
occurrences of $F_1$ occur consecutively.  Repeated applications of
Lemma \ref{Lem:2.2.6} then show that $\mu'_1 = \l_1$.  Lemma
\ref{Lem:2.2.7} (ii) then applies again to yield
$$
E_1^c F_1^c 1_{\mu'} = z'' 1_{\mu'} ,
$$ 
and $z''$ is nonzero because by condition (iv) of Lemma
\ref{Lem:2.4.4}, $c \le \l_1 - 1 < \l_1 = \mu'_1$.  Once again, the
assertion follows by induction in this case.

Finally, we observe that since both $M$ and $\s(M)$ avoid occurrences
of $E_n$ and $F_n$, they lie in the subalgebra of $T$ corresponding to
$\bS_v(n, r)$.  This means that $M$ and $\s(M)$ lie in $\Im(\iota)$,
and the proof follows.
\end{proof}

\begin{cor}\label{Cor:2.4.6}
If $\l \in \Lambda(n, r)$, the elements $E_n 1_\l$, $1_\l F_n$
and $1_\l$ lie in the image of $\iota$.
\end{cor}

\begin{proof}
If $\l_1 = 0$, then $E_n 1_\l = 1_\l F_n = 0$ and the assertion is
trivial.  Otherwise, the assertion follows by combining lemmas 2.4.3
and 2.4.5.
\end{proof}

\begin{proof}[Proof of Theorem \ref{Thm:1.6.1} 
(modulo Lemma \ref{Lem:2.4.4})] By Lemma \ref{Lem:2.4.2} and Corollary
\ref{Cor:2.4.6}, the generators of $T$ listed in Lemma 2.4.1 all lie
in $\Im(\iota)$.  This proves that $\iota$ is surjective, and taken in
conjunction with Proposition \ref{Prop:2.3.11}, we see that $\iota$ is
an isomorphism.  This completes the proof of Theorem 1.6.1 (modulo
Lemma \ref{Lem:2.4.4}).
\end{proof}

\subsection{Proof of Lemma \ref{Lem:2.4.4}}\label{ss:2.5}

The only other ingredient needed to prove Theorem \ref{Thm:1.6.1} is
Lemma 2.4.4.

\begin{defn}\label{Def:2.5.1}
Let $\mu, \nu \in \Lambda(n, r)$.  We say that $\mu$ and $\nu$ are
{\it Z-equivalent} if they become equal after their zero parts have
been deleted.  (In other words, $\mu$ and $\nu$ correspond to the same
parabolic subgroup of the symmetric group.)
\end{defn}

\begin{exam}\label{Ex:2.5.2}
Let $n = 6$, $r = 5$, $\l = (0, 2, 1, 0, 2, 0)$ and $\mu = (2, 0, 1,
2, 0, 0)$.  After deletion of zero parts, each of $\l$ and $\mu$
reduces to $(2, 1, 2)$, so $\l$ and $\mu$ are Z-equivalent.
\end{exam}

\begin{lem}\label{Lem:2.5.3}
Let $\l \in \Lambda(n, r)$ with $\l_1 > 0$.  
\par\noindent{\rm (i)}
There exists a nonzero distinguished monomial $M_1$ in the generators 
$E_2, E_3,$ $\ldots, E_{n-1}$, such that the occurrences of 
$E_{n-1}$ occur consecutively, satisfying 
$$
M_1 = 1_\mu M_1 1_\l.
$$   
Here, $\mu = \mu(\l)$ is such that (a) $\mu$ and $\l$ are Z-equivalent
and (b) $\ell(\mu)$ (see Definition \ref{Def:1.3.1}) is of the form 
$$
(\underbrace{1, 1, \ldots, 1}_{\mu_1},
\underbrace{2, 2, \ldots, 2}_{\mu_2}, \ldots,
\underbrace{k, k, \ldots, k}_{\mu_k}) .
$$  
Furthermore, $\mu_1 = \l_1$, $k \le r$ and thus $\mu_{k+1} = \mu_{k+2}
= \cdots = \mu_{r+1} = \cdots = \mu_n = 0$.
\par\noindent{\rm (ii)} 
Let $\mu$ be as in part (i) above, and let $\nu$ be the unique element of 
$\Lambda(n, r)$ such that 
(a) $\nu_a = \mu_{i+1}$ whenever 
$$
a = 1 + \sum_{j = 1}^i \mu_j
$$ 
for any $0 \le i < r$, and (b) $\nu_a = 0$ for other values of $a$.
(In particular, $\nu_1 = \mu_1 = \l_1$, and $\nu$ and $\mu$ are 
$Z$-equivalent.)  Then there exists a nonzero distinguished 
monomial $M_2$ in the generators $F_2, F_3, \ldots, F_{n-2}$ satisfying 
$$
M_2 = 1_\nu M_2 1_\mu.
$$
\par\noindent{\rm (iii)} 
Let $\mu$ and $\nu$ be as in (ii) above.  Then for each 
$1 \le i \le r$ we have 
$$
\sum_{j = 1}^{\mu_i} \nu_{b(i) + j} = \mu_i = \nu_{b(i) + 1},
$$ 
where $b(i) = \sum_{k < i} \mu_k$.  Summing over all $i$, this yields 
$$
\sum_{j = 1}^r \nu_j = r,
$$ 
and hence $\nu_{r+1} = \nu_{r+2} = \cdots = \nu_n = 0$.
\par\noindent{\rm (iv)} 
Let $\nu$ be as in part (ii) above.  Then there exists a nonzero 
distinguished monomial $M_3$ in the generators $F_1, F_2, F_3, \ldots, 
F_{n-2}$ satisfying 
$$
M_3 = 1_\w M_3 1_\nu.
$$  
Furthermore, the occurrences of $F_1$ occur consecutively, and there
are $\l_1 - 1$ occurrences of $F_1$.
\end{lem}

{\em Note}
Note that in the above situation, $\mu$ and $\nu$ are uniquely determined 
by $\l$.  For example, if $n = 9$, $r = 7$ and 
$$
\l = (2, 0, 0, 3, 0, 0, 0, 0, 2),
$$ then we have 
$$
\mu = (2, 3, 2, 0, 0, 0, 0, 0, 0),
$$ and 
$$ 
\nu = (\underbrace{2, 0}_2, \underbrace{3, 0, 0}_3, 
\underbrace{2,0}_2, 0, 0).
$$  
In this case, we could take 
\begin{align*}
M_1 &= 1_\mu E_3^2 E_4^2 E_5^2 E_6^2 E_7^2 E_8^2 E_2^3 E_3^3 1_\l,\\
M_2 &= 1_\nu F_2^3 F_5^2 F_4^2 F_3^2 1_\mu,\\
M_3 &= 1_\w F_1 F_4 F_3^2 F_6 1_\nu.
\end{align*}

\begin{proof}
To prove (i), let $\l \in \Lambda(n, r)$ be such that $\l_i = 0$ and
$\l_{i+1} = c > 0$.  Lemma \ref{Lem:2.2.7} (i) implies that $E_i^c
1_\l$ is a nonzero element of $T$, and iterated applications of Lemma
\ref{Lem:2.2.6}(i) show that $E_i^c 1_\l = 1_\xi E_i^c$, where
$$
\xi_j = \begin{cases}
 \lambda_{i+1} & \text{ if } j = i,\\
 \lambda_i & \text{ if } j = i+1,\\
 \lambda_j & \text{ otherwise.}
        \end{cases}
$$ 
Repeated applications of this procedure can be used to move all the
zero parts of $\l$ to the right.  Since $\l_1 > 0$ by assumption, we
have $\mu_1 = \l_1$, and it is never necessary to use an application
of $E_1$.  It is necessary to use applications of $E_{n-1}$ if and
only if $\l_n > 0$, but in this case they may all be applied
consecutively.  Since there are at most $r$ nonzero parts in $\l_1$
and since $n > r$, we have $k \le r$.  We have $M_1 \ne 0$ based on an
inductive argument using Lemma \ref{Lem:2.2.11}.  The other assertions
of (i) follow.

The claims in (ii) concerning $\nu$ are routine apart from the
assertion regarding the distinguished monomial.  The proof of (iii)
then follows from the characterization of $\nu$ given in (ii).  (The
reader may find it helpful here to look at the note preceding this
proof.)  The entries $\nu_{r+1}, \nu_{r+2}, \ldots, \nu_n$ can
effectively be ignored for the rest of the proof.

The proof of the last assertion of (ii) follows similar lines to the
proof of (i).  The main difference is that the aim is to move certain
of the zero components of $\l$ to the left.  The basic step involves
$\mu \in \Lambda(n, r)$ such that $\mu_i = c > 0$ and $\mu_{i+1} = 0$,
in which case Lemma \ref{Lem:2.2.7} (ii) shows that $F_i^c 1_\mu$ is a
nonzero element of $T$, and iterated applications of Lemma
\ref{Lem:2.2.6} (ii) show that $F_i^c 1_\mu = 1_\xi F_i^c$, where
$$ 
\xi_j = 
\begin{cases}
 \mu_{i+1} & \text{ if } j = i,\\ 
 \mu_i & \text{ if } j = i+1,\\
 \mu_j & \text{otherwise.}
\end{cases}
$$  
Since $\nu_1 = \mu_1 = \l_1 > 0$ by the definition of $\nu_a$, no
applications of $F_1$ are necessary.  The fact that $\nu_n = 0$ shows
why no applications of $F_{n-1}$ are necessary.  As before, Lemma
2.2.11 shows why $M_2 \ne 0$.

To prove (iv), observe that $\nu$ can be written as the concatenation
of maximal segments of the form
$$
(\nu_i, \underbrace{0, 0, \ldots, 0}_{\nu_i - 1}).
$$  
Let us first deal with the case where $\nu_1 = r$.  Here, the monomial
may be given explicitly as 
$$
M_3 = 1_\w 
F_{\nu_1-1}^1 F_{\nu_1-2}^2 \cdots F_2^{\nu_1 - 2} F_1^{\nu_1 - 1}
1_\nu.
$$ 
This monomial is easily checked to be distinguished.  In the
general case, there is a monomial in the $F_i$ corresponding to each
of the maximal segments mentioned above, and monomials corresponding
to distinct segments commute.  When these monomials are concatenated
and flanked by $1_\w$ and $1_\nu$, we obtain $M_3$.  The occurrences
of $F_1$ all correspond to the segment containing $\nu_1$, and the
explicit formula above shows that these occurrences are consecutive
and that there are $\l_1 - 1$ of them.  (This number may be zero.)  As
in the proof of (ii) above, the fact that $\w_n = 0$ explains why
there are no occurrences of $F_{n-1}$.  Also as above, Lemma
\ref{Lem:2.2.11} shows why $M_3 \ne 0$.
\end{proof}

\begin{proof}[Proof of Lemma \ref{Lem:2.4.4}]
In the notation of Lemma \ref{Lem:2.5.3}, the required monomial is 
$$
M = 1_\w M_3 1_\nu M_2 1_\mu M_1 1_\l.
$$ Properties (i)--(iv) of Lemma \ref{Lem:2.4.4} follow from the
various parts of Lemma \ref{Lem:2.5.3}.  Since $M_1$, $M_2$ and $M_3$
are nonzero, $M$ is nonzero by Lemma \ref{Lem:2.2.11}.

This completes the proof of Lemma \ref{Lem:2.4.4}, and therefore of
Theorem \ref{Thm:1.6.1}.
\end{proof}

\subsection{An alternative presentation of $\svnr$} \label{ss:2.6}

Lusztig \cite[Part IV]{L7} has defined a modified form of a quantized
enveloping algebra, by replacing the zero part of the algebra with an
infinite system of pairwise orthogonal idempotents, acting on modules
as weight space projectors. The modified form has a canonical basis
with remarkable properties, similar to properties of the canonical
basis of the positive part of the original quantized algebra. 

The following presentation of the algebra $\svnr$, compatible with
Lusztig's modified form of $\ugn$, is equivalent with the presentation
given in Theorem \ref{Thm:1.6.1}.

\begin{thm} \label{Thm:2.6.1} 
Assume that $n > r$.  Over $\qv$, the algebra $\svnr$ is isomorphic
with the associative algebra (with 1) given by the generators $\ii_\l$
($\lambda \in \Lambda(n,r)$), $\EE_i, \FF_i$ ($1 \le i \le n$) and
relations
\begin{align}
\ii_\l \, \ii_\m &= \delta_{\l,\m} \, \ii_\l; \qquad \sum_{\l \in
  \Lambda(n,r)} \ii_\l = 1; \tag{R1}\\
\EE_i\, \ii_\l &= \begin{cases}
\ii_{\lambda + \a_i} \EE_i & \text{ if } \lambda_{i+1} > 0, \tag{R2}\\
0 & \text{ otherwise;}
\end{cases} \\
\FF_i \,\ii_\l &= \begin{cases}
 \ii_{\lambda - \a_i} \FF_i & \text{ if } \lambda_{i} > 0, \tag{R3} \\
 0 & \text{ otherwise;}
\end{cases}\\
\EE_i \FF_j - \FF_j \EE_i &= \delta_{i,j} \sum_{\l \in \Lambda(n,r)}
[\l_j - \l_{j+1}]\,\ii_\l \tag{R4}
\end{align} 
along with relations (Q6)--(Q9) of Definition \ref{Def:1.4.1}. Here we
regard weights as infinite periodic sequences, as in Section
\ref{ss:2.2} above.  
\end{thm}

\begin{proof}
Let $A$ be the algebra defined by the presentation of the theorem, and
let $T$ be the algebra defined by the presentation of Theorem
\ref{Thm:1.6.1}.

By Definition \ref{Def:2.2.2}, Lemma \ref{Lem:2.2.4}, and Lemma
\ref{Lem:2.2.6} the elements $1_\l$ ($\l \in \Lambda(n,r)$), $E_i,
F_i$ ($1\le i \le n$) of $T$ satisfy the relations (R1)--(R3).  The
$E_i, F_i$ ($1\le i \le n$) already satisfy relations (Q6)--(Q9) by
assumption. By applying Remark \ref{Rem:1.6.4} to the results in
\cite[Theorem 3.4]{DG} we see that (R4) holds as well for all $1 \le
i,j < n$. If one or both of $i,j$ is equal to $n$, then we choose a
different embedding of $\bS_v(n,r)$ in $T$, one which includes the
values of $i,j$ in question, and again apply Remark \ref{Rem:1.6.4} to
the results in \cite[Theorem 3.4]{DG} to see that (R4) holds in that
case as well. In $T$ we have by Lemma \ref{Lem:2.2.4} the equalities
\[
 K_i = K_i \sum_\l 1_\l = \sum_\l v^{\l_i} 1_\l; \quad K_i^{-1} =
 K_i^{-1} \sum_\l 1_\l = \sum_\l v^{-\l_i} 1_\l  \tag{\text{$\ast$}}
\]
for any $1 \le i \le n$, where the sums are taken over all $\lambda
\in \Lambda(n,r)$.  Hence, the elements $1_\l$ ($\l \in
\Lambda(n,r)$), $E_i, F_i$ ($1\le i \le n$) generate $T$, and the
map
\[
  \ii_\l \to 1_\l, \quad \EE_i \to E_i, \quad \FF_i \to F_i 
\]
defines a surjective quotient map from $A$ onto $T$.

On the other hand, in the algebra $A$ one defines elements $\KK_i =
\sum_\l v^{\l_i} \ii_\l$, $\KK_i^{-1} = \sum_\l v^{-\l_i} \ii_\l$.  By
following the same line of argument as in the proof of \cite[Theorem
3.4]{DG}, these elements, along with the elements $\EE_i$, $\FF_i$ for
$1 \le i < n$, satisfy the defining relations (Q1)--(Q9), (Q15), (Q16)
of $T$. It remains to show that the elements $\EE_n,\FF_n$ also
satisfy those relations (along with the $\KK_i, \KK_i^{-1}$). Only
relations (Q3), (Q4), and (Q5) are in question since the other
relations either hold by assumption or do not involve the elements
$\EE_n,\FF_n$.

We now verify that relation (Q3) holds for $\EE_n$. By definition of
$\KK_i$ we have
\[
\KK_i \EE_n = \sum_\l v^{\l_i} \, \ii_\l \, \EE_n
\]
and by relation (R2) this takes the form
\[
\KK_i \EE_n = \sum_\l v^{\l_i} \, \EE_n \, \ii_{\l-\a_n}
\]
where, for convenience of notation, we take both sums over the set of
all $\l \in \Z^n$ satisfying $\sum \l_i = r$, with the understanding
that $\ii_\l$ is interpreted to be 0 in case any part of $\l$ is
negative. (This makes all the sums in question finite.) Now replace
$\l - \a_n$ by $\mu \in \Z^n$ and the above gives relation (Q3) for
$\EE_n$.

The proof that relation (Q4) holds for $\FF_n$ is similar. 

Finally, we verify that (Q5) holds.  By the given relation (R4) we have
\[
 \EE_i \FF_j - \FF_j \EE_i = \delta_{ij} \sum_\l [\l_j - \l_{j+1}]\, \ii_\l
\]
and this gives 0 unless $i=j$, so (Q5) holds in case $i\ne
j$. Assuming that $i=j$,  the above sum becomes 
\begin{align*}
\EE_i \FF_i - \FF_i \EE_i &= \sum_\l 
 \frac{v^{\l_i - \l_{i+1}} - v^{-\l_i + \l_{i+1}}}{v-v^{-1}} \,\ii_\l \\
&= \frac{(\sum_\l v^{\l_i}\ii_\l v^{-\l_{i+1}}\ii_\l) - 
   (\sum_\l v^{-\l_i}\ii_\l v^{\l_{i+1}}\ii_\l)}{v-v^{-1}}  \\
&= \frac{(\sum_\l v^{\l_i}\ii_\l)(\sum_\l v^{-\l_{i+1}}\ii_\l) - 
   (\sum_\l v^{-\l_i}\ii_\l)(\sum_\l v^{\l_{i+1}}\ii_\l)}{v-v^{-1}} 
\end{align*}
using the orthogonality of the system of idempotents. By the
definition of the $\KK_i, \KK_i^{-1}$ this proves (Q5) in case $i=j$.

We claim that the elements $\KK_i, \KK_i^{-1}$, $\EE_i, \FF_i$ (for $1
\le i \le n$) generate $A$.  To see this, it suffices to show that the
$\KK_i, \KK_i^{-1}$ generate the zero part of $A$ (the span of the
$\ii_\l$). From the definition of $\KK_i$ and $\KK_i^{-1}$ it follows
that
\[
 \KK_j\, \ii_\l = v^{\l_j}\, \ii_\l; \quad 
 \KK_j^{-1}\, \ii_\l = v^{-\l_j}\, \ii_\l
\]
and thus $\KK_j = \KK_j \sum_\l \ii_\l = \sum_\l v^{\l_j} \, \ii_\l$
and $\KK_j^{-1} = \KK_j^{-1}\sum_\l \ii_\l = \sum_\l v^{-\l_j} \,
\ii_\l$, where the sums are over all $\l \in \Lambda(n,r)$. Hence it
follows that
\begin{align*}
  \sqbinom{\KK_j}{t} &= \prod_{s=1}^t  
    \frac{\KK_j v^{-s+1} - \KK_j^{-1} v^{s-1}}{v^s - v^{-s}} \\
 &= \prod_{s=1}^t  
    \frac{ (\sum_\l v^{\l_j-s+1}\,\ii_\l) - (\sum_\l
    v^{-\l_j+s-1}\,\ii_\l) } {v^s - v^{-s}} \\
&= \prod_{s=1}^t \sum_\l \frac{ v^{\l_j-s+1} - v^{-\l_j+s-1}}
    {v^s-v^{-s}} \; \ii_\l \\
&= \sum_\l \prod_{s=1}^t \frac{ v^{\l_j-s+1} - v^{-\l_j+s-1}}
    {v^s-v^{-s}} \; \ii_\l \\
&= \sum_\l \sqbinom{\l_j}{t} \, \ii_\l
\end{align*} 
where we have again made use of the orthogonality of the idempotents
to interchange the product and sum. From this and the othogonality of
idempotents it follows that for any $\mu \in \Lambda(n,r)$ we have
\begin{align*}
\sqbinom{\KK_1}{\mu_1}\cdots \sqbinom{\KK_n}{\mu_n} &= 
\prod_{j=1}^n \left( \sum_\l \sqbinom{\l_j}{\mu_j}\, \ii_\l \right) \\
&= 
\sum_\l \sqbinom{\l_1}{\mu_1}\cdots \sqbinom{\l_n}{\mu_n}  \, \ii_\l 
\end{align*}
where $\l$ runs over the set $\Lambda(n,r)$ in the sums. The only
non-zero term in the last sum is when $\l = \mu$, so 
\[
  \sqbinom{\KK_1}{\mu_1}\cdots \sqbinom{\KK_n}{\mu_n} = \ii_\mu.
\]
This proves the claim. (The reader should refer to \cite[\S1.3]{L7}
for definitions and basic properties of quantized binomial
coefficients used here.)

We have shown that the elements $\KK_i, \KK_i^{-1}$, $\EE_i$, $\FF_i$ 
(for $1 \le i \le n$) generate the algebra $A$, and moreover satisfy
all the defining relations for the algebra $T$. 
It follows that the map
\[
  K_i^{\pm 1} \to \KK_i^{\pm 1},\quad E_i \to \EE_i, \quad F_i
  \to \FF_i
\]
is a surjective quotient map from $T$ onto $A$. 

Now consider the composite surjective map $T \to A \to T$. This is the
clearly identity on $E_i, F_i$.  Moreover, by equations ($\ast$) above
the composite map takes $K_i$ to $\sum_\l v^{\l_i}1_\l =
K_i$. Similarly, it takes $K_i^{-1}$ to itself.  Thus the composite is
identity, and thus each quotient map $T\to A$ and $A \to T$ is an
algebra isomorphism.
\end{proof}

\section{The classical case}
All of the results of this paper have analogues in the case $v=1$.
The proofs run parallel to the arguments given here, but are often
easier. We will outline the main results here, and leave it to the
reader to fill in the details.

\subsection{The affine Schur algebra} \label{ss:3.1} 
The analogue of Definition \ref{Def:1.3.4} is the following

\begin{defn}\label{Def:3.1.1} 
The affine Schur algebra $\shatnr$ over $\Z$ is defined by
\[
\shatnr := \End_\W \left( \bigoplus_{\l \in \Lambda(n, r)} x_\l \W
\right),
\]
where $x_\l = \sum_{w \in \W_\l} w$. 
\end{defn}

There is a basis of $\shatnr$ similar to the basis of $\sqnr$ given in
Theorem \ref{Thm:1.3.6}. The details are left to the reader.

\begin{defn} \label{Def:3.1.2}
The associative, unital algebra $U(\widehat{\gl_n})$ over $\mathbb{Q}$
is given by generators
\[
  e_i, f_i, H_i  \quad (1 \le i \le n)
\]
subject to the relations
\begin{align}
H_iH_j &= H_jH_i; \tag{q1}\\
H_ie_j-e_jH_i &= \epsilon^+(i,j) e_j; \tag{q2}\\
H_if_j-f_jH_i &= \epsilon^-(i,j) f_j; \tag{q3}\\
e_if_j-f_je_i &= \delta_{ij} (H_j-H_{j+1}); \tag{q4}\\
e_ie_j &= e_je_i \text{ if $i$ and $j$ are not adjacent}; \tag{q5}\\
f_if_j &= f_jf_i \text{ if $i$ and $j$ are not adjacent}; \tag{q6}\\
e_i^2e_j-2e_ie_je_i+e_je_i^2 &= 0 \text{ if $i$ and $j$ are
  adjacent}; \tag{q7}\\
f_i^2f_j-2f_if_jf_i+f_jf_i^2 &= 0 \text{ if $i$ and $j$ are
  adjacent}. \tag{q8}
\end{align}
As in Definition \ref{Def:1.4.2} the notion of adjacency takes place
in the Dynkin diagram of type $\widehat{A}_{n-1}$, so we read indices
modulo $n$ in this definition.
\end{defn}

This algebra is a Hopf algebra in a natural way, and the quotient of
$U(\widehat{\gl_n})$ by the kernel of its action on a suitably defined
tensor space is isomorphic as a $\mathbb{Q}$-algebra to
$\mathbb{Q}\otimes_\Z \shatnr$.

\subsection{Main results}\label{ss:3.2}

The analogue of Theorem \ref{Thm:1.6.1} is the
\begin{thm}\label{Thm:3.2.1}
Let $n > r$, and identify $\shatnr$ with the quotient of
$U(\widehat{\gl_n})$ acting on tensor space.  Over $\mathbb{Q}$, the
affine Schur algebra $\shatnr$ is given by generators $e_i, f_i,
H_i$ ($1 \le i \le n$) subject to relations (q1) to (q8) of
Definition \ref{Def:3.1.2} (reading indices modulo $n$), together
with the relations
\begin{align}
H_1 + \cdots + H_n &= r;  \tag{q9}\\
H_i(H_i - 1)(H_i - 2)\cdots(H_i - r) &= 0. \tag{q10}
\end{align}
\end{thm}

There is also an equivalent version in terms of idempotents, analogous
to Theorem \ref{Thm:2.6.1}, which we now state.

\begin{thm} \label{Thm:3.2.2} 
Assume that $n > r$.  Over $\mathbb{Q}$, the algebra $\shatnr$ is
isomorphic with the associative algebra (with 1) given by the
generators $\iii_\l$ ($\lambda \in \Lambda(n,r)$), $\eee_i, \fff_i$
($1 \le i \le n$) and relations
\begin{align}
\iii_\l \, \iii_\m &= \delta_{\l,\m} \, \iii_\l; \qquad \sum_{\l \in
  \Lambda(n,r)} \iii_\l = 1; \tag{r1}\\
\eee_i\, \iii_\l &= \begin{cases}
\iii_{\lambda + \a_i} \eee_i & \text{ if } \lambda_{i+1} > 0, \tag{r2}\\
0 & \text{ otherwise;}
\end{cases} \\
\fff_i \,\iii_\l &= \begin{cases}
 \iii_{\lambda - \a_i} \fff_i & \text{ if } \lambda_{i} > 0, \tag{r3} \\
 0 & \text{ otherwise;}
\end{cases}\\
\eee_i \fff_j - \fff_j \eee_i &= \delta_{i,j} \sum_{\l \in \Lambda(n,r)}
(\l_j - \l_{j+1})\,\iii_\l \tag{r4}
\end{align} 
along with relations (q5)--(q8) of Definition \ref{Def:3.1.2}. Here we
regard weights as infinite periodic sequences, as in Section
\ref{ss:2.2} above.
\end{thm}
These relations are obtained from those in Theorem \ref{Thm:2.6.1}
by setting $v=1$.

\bibliographystyle{amsalpha}

\end{document}